\documentclass[11pt, reqno]{amsart}

\usepackage[all]{xy}
\usepackage{latexsym}
\usepackage{amsmath}
\usepackage{mathtools}
\usepackage{amsthm}
\usepackage{amssymb}
\usepackage{a4wide}
\usepackage[utf8]{inputenc}
\usepackage{dsfont}
\usepackage{multirow}
\usepackage{hyperref}
\usepackage{stmaryrd}

\theoremstyle{plain}
\newtheorem{thm}{Theorem}[section]
\newtheorem*{thmintro}{Theorem}
\newtheorem{prop}[thm]{Proposition}
\newtheorem{propA}[thm]{Proposition}
\newtheorem{lem}[thm]{Lemma}
\newtheorem{cor}[thm]{Corollary}

\theoremstyle{definition}
\newtheorem{dei}[thm]{Definition}

\theoremstyle{remark}
\newtheorem*{rem}{\sc Remark}
\newtheorem*{rems}{\sc Remarks}
\newtheorem*{pf}{\sc Proof}

\newcommand{\Rb}{\mathbb{R}}
\newcommand{\Cb}{\mathbb{C}}
\newcommand{\Kb}{\mathbb{K}}
\newcommand{\Sb}{\mathbb{S}}

\newcommand{\Nb}{\mathbb{N}}
\newcommand{\Zb}{\mathbb{Z}}

\newcommand{\A}{\mathcal{A}}
\newcommand{\B}{\mathcal{B}}
\newcommand{\C}{\mathcal{C}}
\newcommand{\D}{\mathcal{D}}

\newcommand{\M}{\mathcal{M}}
\newcommand{\N}{\mathcal{N}}
\renewcommand{\O}{\mathcal{O}}
\renewcommand{\P}{\mathcal{P}}

\newcommand{\R}{\mathcal{R}}

\newcommand{\T}{\mathcal{T}}

\newcommand{\V}{\mathcal{V}}
\newcommand{\W}{\mathcal{W}}
\newcommand{\jo}{\mathring{\jmath}}
\newcommand{\fo}{\mathring{f}}

\newcommand{\Cx}{\mathrm{Cx}}
\newcommand{\Bm}{\mathrm{B}}
\renewcommand{\Im}{\mathrm{I}}
\newcommand{\GL}{\mathrm{GL}}
\newcommand{\Mrm}{\mathrm{M}}
\newcommand{\coor}{\mathrm{coor}}
\newcommand{\Map}{\mathrm{Map}}
\newcommand{\Hm}{\mathrm{H}}

\newcommand{\Vect}{\mathsf{Vect}}

\newcommand{\qCx}{\mathrm{q}\mathrm{Cx}}
\newcommand{\qCb}{\mathrm{q}\underline{\mathbb{C}}}

\newcommand{\cqfd}{\ \hfill \square}

\newcommand{\As}{\mathrm{As}}

\newcommand{\Lie}{\mathrm{Lie}}
\newcommand{\Com}{\mathcal{C}om}

\newcommand{\ash}{\textrm{!`}}

\newcommand{\Mco}{M^{\textrm{coor}}}
\newcommand{\Maf}{M^{\textrm{aff}}}
\newcommand{\tE}{\widetilde{E}}
\newcommand{\bE}{\overline{E}}

\newcommand{\Id}{\operatorname{Id}}
\newcommand{\id}{\operatorname{id}}
\newcommand{\End}{\operatorname{End}}

\newcommand{\Hom}{\operatorname{Hom}}
\newcommand{\sgn}{\operatorname{sgn}}

\newcommand{\jn}{{\tiny \vcenter{\xymatrix@M=0pt@R=6pt@C=4pt{& 1 &\\ & \ar@{-}[dd] &\\ && \\ & \ar@{{*}}[u] &}}}}

\newcommand{\jnr}{{\tiny \vcenter{\xymatrix@M=0pt@R=4pt@C=4pt{& 1 &\\ & \ar@{-}[dddd] &\\ && \\ & \ar@{{*}}[u] \ar@{{*}}[d] & \\ && \\ &&}} + \vcenter{\xymatrix@M=0pt@R=4pt@C=4pt{& 1 &\\ & \ar@{-}[dddd] &\\ && \\ && \\ && \\ &&}}}}

\newcommand{\qjnr}{{\tiny \vcenter{\xymatrix@M=0pt@R=4pt@C=4pt{& 1 &\\ & \ar@{-}[dddd] &\\ && \\ & \ar@{{*}}[u] \ar@{{*}}[d] & \\ && \\ &&}}}}

\newcommand{\jnsi}{\vcenter{\xymatrix@M=0pt@R=6pt@C=4pt{& \ar@{-}[dd] &\\ && \\ & \ar@{{*}}[u] &}}}

\newcommand{\jnsip}{\vcenter{\xymatrix@M=0pt@R=6pt@C=4pt{& \ar@{-}[dd] &\\ && \\ & \ar@{{*}}[u]_' &}}}

\newcommand{\jnrsi}{\vcenter{\xymatrix@M=0pt@R=4pt@C=4pt{& \ar@{-}[dddd] &\\ && \\ & \ar@{{*}}[u] \ar@{{*}}[d] & \\ && \\ &&}} + \vcenter{\xymatrix@M=0pt@R=4pt@C=4pt{& \ar@{-}[dddd] &\\ && \\ && \\ && \\ &&}}}

\newcommand{\qjnrsi}{\vcenter{\xymatrix@M=0pt@R=4pt@C=4pt{& \ar@{-}[dddd] &\\ && \\ & \ar@{{*}}[u] \ar@{{*}}[d] & \\ && \\ &&}}}

\newcommand{\lie}{{\tiny \vcenter{\xymatrix@M=0pt@R=4pt@C=4pt{1 && 2\\ && \\ & *{} \ar@{-}[lu] \ar@{-}[d] \ar@{-}[ru] & \\ &&}}}}

\newcommand{\jac}{{\tiny \vcenter{\xymatrix@M=0pt@R=4pt@C=4pt{1 && 2 && 3\\ &&&& \\ & \ar@{-}[lu] \ar@{-}[dr] \ar@{-}[ru] &&& \\ && *{} \ar@{-}[d] \ar@{-}[urur] &&\\ &&&&}} + \vcenter{\xymatrix@M=0pt@R=4pt@C=4pt{2 && 3 && 1\\ &&&& \\ & \ar@{-}[lu] \ar@{-}[dr] \ar@{-}[ru] &&& \\ && *{} \ar@{-}[d] \ar@{-}[urur] &&\\ &&&&}} + \vcenter{\xymatrix@M=0pt@R=4pt@C=4pt{3 && 1 && 2\\ &&&& \\ & \ar@{-}[lu] \ar@{-}[dr] \ar@{-}[ru] &&& \\ && *{} \ar@{-}[d] \ar@{-}[urur] &&\\ &&&&}}}}

\newcommand{\jlr}{{\tiny \vcenter{\xymatrix@M=0pt@R=4pt@C=4pt{1 && 2\\ \ar@{-}[dd] &&\\ &&\\ \ar@{{*}}[u] && \\ & *{} \ar@{-}[ul] \ar@{-}[d] \ar@{-}[ur] &\\ &&}}}}

\newcommand{\ljr}{{\tiny \vcenter{\xymatrix@M=0pt@R=4pt@C=4pt{1 && 2\\ \ar@{-}[dr] && \ar@{-}[ld]\\ & *{} \ar@{-}[dd] &\\ &&\\ & \ar@{{*}}[u] &}}}}

\newcommand{\jlrt}{{\tiny \vcenter{\xymatrix@M=0pt@R=4pt@C=4pt{2 && 1\\ \ar@{-}[dd] &&\\ &&\\ \ar@{{*}}[u] && \\ & *{} \ar@{-}[ul] \ar@{-}[d] \ar@{-}[ur] &\\ &&}}}}

\newcommand{\ljrt}{{\tiny \vcenter{\xymatrix@M=0pt@R=4pt@C=4pt{2 && 1\\ \ar@{-}[dr] && \ar@{-}[ld]\\ & *{} \ar@{-}[dd] &\\ &&\\ & \ar@{{*}}[u] &}}}}

\newcommand{\slie}{{\tiny \vcenter{\xymatrix@M=0pt@R=4pt@C=4pt{1 &&& 2\\ &&& \\ && *{} \ar@{-}[llu] \ar@{-}[d] \ar@{-}[ru] & s^{-1} \\ &&}}}}

\newcommand{\sjlr}{{\tiny \vcenter{\xymatrix@M=0pt@R=4pt@C=4pt{1 &&& 2\\ \ar@{-}[dd] &&&\\ &&&\\ \ar@{{*}}[u] &&& \\ && *{} \ar@{-}[ull] \ar@{-}[d] \ar@{-}[ur] & s^{-1}\\ &&&}}}}

\newcommand{\sljr}{{\tiny \vcenter{\xymatrix@M=0pt@R=4pt@C=4pt{1 &&& 2\\ \ar@{-}[drr] &&& \ar@{-}[ld]\\ && *{} \ar@{-}[dd] & s^{-1} \\ &&&\\ && \ar@{{*}}[u] &}}}}

\newcommand{\sljrp}{{\tiny \vcenter{\xymatrix@M=0pt@R=4pt@C=4pt{1 &&& 2\\ \ar@{-}[drr] &&& \ar@{-}[ld]\\ && *{} \ar@{-}[dd] & s^{-1} \\ &&&\\ && \ar@{{*}}[u]_' &}}}}

\newcommand{\sjac}{{\tiny \vcenter{\xymatrix@M=0pt@R=4pt@C=4pt{1 && 2 && 3\\ &&&& \\ s^{-1} & \ar@{-}[lu] \ar@{-}[dr] \ar@{-}[ru] &&& \\ & s^{-1} & *{} \ar@{-}[d] \ar@{-}[urur] & \ &\\ &&&&}} + \vcenter{\xymatrix@M=0pt@R=4pt@C=4pt{2 && 3 && 1\\ &&&& \\ s^{-1} & \ar@{-}[lu] \ar@{-}[dr] \ar@{-}[ru] &&& \\ & s^{-1} & *{} \ar@{-}[d] \ar@{-}[urur] & \ &\\ &&&&}} + \vcenter{\xymatrix@M=0pt@R=4pt@C=4pt{3 && 1 && 2\\ &&&& \\ s^{-1} & \ar@{-}[lu] \ar@{-}[dr] \ar@{-}[ru] &&& \\ & s^{-1} & *{} \ar@{-}[d] \ar@{-}[urur] & \ &\\ &&&&}}}}

\newcommand{\sjlrt}{{\tiny \vcenter{\xymatrix@M=0pt@R=4pt@C=4pt{2 &&& 1\\ \ar@{-}[dd] &&&\\ &&&\\ \ar@{{*}}[u] &&& \\ && *{} \ar@{-}[ull] \ar@{-}[d] \ar@{-}[ur] & s^{-1}\\ &&&}}}}

\newcommand{\sljrt}{{\tiny \vcenter{\xymatrix@M=0pt@R=4pt@C=4pt{2 &&& 1\\ \ar@{-}[drr] &&& \ar@{-}[ld]\\ && *{} \ar@{-}[dd] & s^{-1} \\ &&&\\ && \ar@{{*}}[u] &}}}}

\newcommand{\as}{{\tiny \vcenter{\xymatrix@M=2pt@R=2pt@C=4pt{&& \\ & \overline{\mu} \ar@{-}[lu] \ar@{-}[d] \ar@{-}[ru] & \\ &&}}}}

\newcommand{\ass}{{\tiny \vcenter{\xymatrix@M=2pt@R=2pt@C=2pt{&&&& \\ & \overline{\mu} \ar@{-}[lu] \ar@{-}[dr] \ar@{-}[ru] &&& \\ && \overline{\mu} \ar@{-}[d] \ar@{-}[urur] & \ \ &\\ &&&&}} + \vcenter{\xymatrix@M=2pt@R=2pt@C=2pt{&&&& \\ &&& \overline{\mu} \ar@{-}[lu] \ar@{-}[dl] \ar@{-}[ru] & \\ & \ \ & \overline{\mu} \ar@{-}[d] \ar@{-}[ulul] &&\\ &&&&}}}}

\newcommand{\sojn}{{\tiny \vcenter{\xymatrix@M=0pt@R=6pt@C=4pt{& 1 &\\ & \ar@{-}[d] &\\ & \circ & \\ & \ar@{-}[u] &}}}}

\newcommand{\soqjnr}{{\tiny \vcenter{\xymatrix@M=0pt@R=4pt@C=4pt{& 1 &\\ & \ar@{-}[d] &\\ & \circ & \\ & *{} \ar@{-}[u] \ar@{-}[d] & \\ & \circ \ar@{-}[d] & \\ &&}}}}

\newcommand{\sojlr}{{\tiny \vcenter{\xymatrix@M=0pt@R=4pt@C=4pt{1 && 2\\ \ar@{-}[d] &&\\ \circ &&\\ \ar@{-}[u] && \\ & *{} \ar@{-}[ul] \ar@{-}[d] \ar@{-}[ur] &\\ &&}}}}

\newcommand{\soljr}{{\tiny \vcenter{\xymatrix@M=0pt@R=4pt@C=4pt{1 && 2\\ \ar@{-}[dr] && \ar@{-}[ld]\\ & *{} \ar@{-}[d] &\\ & \circ &\\ & \ar@{-}[u] &}}}}

\newcommand{\jjac}{{\tiny \vcenter{\xymatrix@M=0pt@R=4pt@C=4pt{1 && 2 && 3\\ &&&&\\ \circ \ar@{-}[u] &&&&\\ \ar@{-}[u] &&&&\\ & \ar@{-}[lu] \ar@{-}[dr] \ar@{-}[ru] &&&\\ && *{} \ar@{-}[d] \ar@{-}[urur] &&\\ &&&&}} + \vcenter{\xymatrix@M=0pt@R=4pt@C=4pt{2 && 3 && 1\\ &&&&\\ &&&& \circ \ar@{-}[u] \\ &&&& \ar@{-}[u]\\ & \ar@{-}[lu] \ar@{-}[dr] \ar@{-}[ru] &&&\\ && *{} \ar@{-}[d] \ar@{-}[urur] &&\\ &&&&}} + \vcenter{\xymatrix@M=0pt@R=4pt@C=4pt{3 && 1 && 2\\ &&&&\\ && \circ \ar@{-}[u] &&\\ && \ar@{-}[u] &&\\ & \ar@{-}[lu] \ar@{-}[dr] \ar@{-}[ru] &&&\\ && *{} \ar@{-}[d] \ar@{-}[urur] &&\\ &&&&}} + \vcenter{\xymatrix@M=0pt@R=4pt@C=4pt{1 && 2 & 3 &\\ &&&& \\ & *{} \ar@{-}[lu] \ar@{-}[d] \ar@{-}[ru] &&&\\ & \circ \ar@{-}[d] &&& \\ & \ar@{-}[dr] &&&\\ && *{} \ar@{-}[d] \ar@{-}[ur] &&\\ &&&&}} + \vcenter{\xymatrix@M=0pt@R=4pt@C=4pt{3 && 1 & 2 &\\ &&&& \\ & *{} \ar@{-}[lu] \ar@{-}[d] \ar@{-}[ru] &&&\\ & \circ \ar@{-}[d] &&& \\ & \ar@{-}[dr] &&&\\ && *{} \ar@{-}[d] \ar@{-}[ur] &&\\ &&&&}} + \vcenter{\xymatrix@M=0pt@R=4pt@C=4pt{1 && 2 && 3\\ &&&& \\ & \ar@{-}[lu] \ar@{-}[dr] \ar@{-}[ru] &&& \\ && *{} \ar@{-}[d] \ar@{-}[urur] &&\\ && \circ \ar@{-}[d] &&\\ &&&&}} + \vcenter{\xymatrix@M=0pt@R=4pt@C=4pt{2 && 3 && 1\\ &&&& \\ & \ar@{-}[lu] \ar@{-}[dr] \ar@{-}[ru] &&& \\ && *{} \ar@{-}[d] \ar@{-}[urur] &&\\ && \circ \ar@{-}[d] &&\\ &&&&}} + \vcenter{\xymatrix@M=0pt@R=4pt@C=4pt{3 && 1 && 2\\ &&&& \\ & \ar@{-}[lu] \ar@{-}[dr] \ar@{-}[ru] &&& \\ && *{} \ar@{-}[d] \ar@{-}[urur] &&\\ && \circ \ar@{-}[d] &&\\ &&&&}}}}

\newcommand{\jjlr}{{\tiny \vcenter{\xymatrix@M=0pt@R=4pt@C=4pt{1 && 2\\ \ar@{-}[d] &&\\ \circ &&\\ \circ \ar@{-}[u] \ar@{-}[d] && \\ &&\\ & *{} \ar@{-}[ul] \ar@{-}[d] \ar@{-}[ur] &\\ &&}} + \vcenter{\xymatrix@M=0pt@R=4pt@C=4pt{1 && 2\\ \ar@{-}[d] &&\\ \circ &&\\ \ar@{-}[u] && \\ & *{} \ar@{-}[ul] \ar@{-}[d] \ar@{-}[ur] &\\ & \circ \ar@{-}[d] &\\ &&}} + \vcenter{\xymatrix@M=0pt@R=4pt@C=4pt{1 && 2\\ \ar@{-}[dr] && \ar@{-}[ld]\\ & *{} \ar@{-}[d] &\\ & \circ &\\ & \circ \ar@{-}[u] \ar@{-}[d] &\\ &&}}}}

\newcommand{\jljr}{{\tiny \vcenter{\xymatrix@M=0pt@R=4pt@C=4pt{1 && 2\\ \ar@{-}[d] && \ar@{-}[d] \\ \circ \ar@{-}[d] && \circ \ar@{-}[d] \\ &&\\ & *{} \ar@{-}[ul] \ar@{-}[d] \ar@{-}[ur] &\\ &&}} + \vcenter{\xymatrix@M=0pt@R=4pt@C=4pt{1 && 2\\ && \ar@{-}[d]\\ && \circ \\ && \ar@{-}[u] \\ & *{} \ar@{-}[ul] \ar@{-}[d] \ar@{-}[ur] &\\ & \circ \ar@{-}[d] &\\ &&}} - \vcenter{\xymatrix@M=0pt@R=4pt@C=4pt{1 && 2\\ \ar@{-}[d] &&\\ \circ &&\\ \ar@{-}[u] && \\ & *{} \ar@{-}[ul] \ar@{-}[d] \ar@{-}[ur] &\\ & \circ \ar@{-}[d] &\\ &&}}}}

\newcommand{\sojnsi}{{\tiny \vcenter{\xymatrix@M=0pt@R=4pt@C=4pt{& \ar@{-}[d] &\\ & \circ & \\ & \ar@{-}[u] &}}}}

\newcommand{\soqjnrsi}{{\tiny \vcenter{\xymatrix@M=0pt@R=4pt@C=4pt{& \ar@{-}[d] &\\ & \circ & \\ & *{} \ar@{-}[u] \ar@{-}[d] & \\ & \circ \ar@{-}[d] & \\ &&}}}}


\newcommand{\draftnote}[1]{}

\author{Joan Mill\`es}

\title{Complex manifolds as families of homotopy algebras}

\begin{document}

\begin{abstract}
We prove an equivalence of categories from formal complex structures with formal holomorphic maps to homotopy algebras over a simple operad with its associated homotopy morphisms. We extend this equivalence to complex manifolds. A complex structure on a smooth manifold corresponds in this way to a family of algebras indexed by the points of the manifold.
\end{abstract}

\maketitle

\tableofcontents

\section*{Introduction}

A \emph{complex manifold} is a manifold $M$ endowed with an atlas of charts to $\Cb^n$ for a certain fixed $n \in \Nb$, such that the transition maps are holomorphic. By a theorem of A. Newlander and L. Nirenberg \cite{NewlanderNirenberg}, such a structure is described equivalently as an almost complex structure $J : TM \rightarrow TM$ satisfying an integrability condition. We use this equivalent formulation to propose a third description: a complex structure is a family of algebras over a certain operad $\Cx_{\infty}$.\\

We make use of the notion of operad to encode algebras: the latter are representations of a given operad. Examples of operads are associative algebras and, in this case, representations are modules over the given associative algebra. For instance, representations of the algebra of complex numbers $\Cb$ are $\Cb$-vector spaces. Another example, which is not an associative algebra, is the operad $\Lie$ whose representations are Lie-algebras.\\

This new description of complex structures makes possible the study of problems in complex geometry with the tools of homological algebra. For example, we get directly notions of cohomology theory and of deformation theory for formal complex structures. Furthermore, algebras fit into the more general context of differential graded (dg) objects and we therefore obtain a notion of dg formal complex manifold. We plan to use this language to describe complex structures with singularities in a future work.\\

In \cite{Merkulov04, Merkulov05, Merkulov06}, S. Merkulov began the description of several geometries in the context of homological algebra (Hertling-Manin, Nijenhuis and Poisson structures). Later, H. Strohmayer dealt with the case of bi-Hamiltonian structures \cite{Strohmayer10}. These notions consist of an underlying smooth manifold endowed with a particular structure. It is sometimes possible to study similarly objects with a geometric flavor such as quantum BV structures \cite{Schwarz93, Khudaverdian04, Merkulov10}. For the aforementioned geometries, the extra structures are described by local rules. In all these cases, we can restrict the local rules defining the extra data to the formal neighborhood of a point \cite{Grothendieck, Kontsevich03}. Since they are described by smooth applications satisfying differential equations, it corresponds to replacing the applications by their infinite Taylor series and the differential equations by the associated algebraic equations between the Taylor coefficients.\\

A powerful tool of homological algebra is the Homotopy Transfer Theorem \cite{Kadeishvili, LodayVallette}. Remarkably, when applied to the algebras corresponding to the quantum BV structures, we obtain precisely the Feynman diagrams appearing in the Batalin-Vilkovisky quantization \cite{Mnev08a, Mnev08b, Merkulov10}. A second illustration of the fruitfulness of the algebraic approach is the reformulation in \cite{Merkulov08} of the universal quantization of Poisson structures in terms of a morphism of props. This approach shows the importance of traces in that context.\\

In \cite{Merkulov05}, S. Merkulov studies \emph{Nijenhuis structures}, that is endomorphisms $J$ on the tangent bundle of a smooth manifold which satisfy an integrability condition. S. Merkulov provides in his article an operad $\mathrm{N}_{\infty}$ whose algebras are formal Nijenhuis structures. H. Strohmayer proves in \cite{Strohmayer} that the operad $\mathrm{N}$, on which the operad $\mathrm{N}_{\infty}$ is based, is a Koszul operad. This shows that $\mathrm{N}_{\infty}$-algebras are homotopy $\mathrm{N}$-algebras and that the two notions of algebras encode the same homotopy categories of algebras. The main example of a Nijenhuis structure is a complex structure, which satisfies in addition the equality $J^{2} = - \Id$ ($J$ is an almost complex structure). It is natural to wonder if formal complex structures can also be modeled similarly. The first goal of this paper is to answer this question.\\

We built a Koszul operad $\Cx$ based on the algebra of complex numbers $\Cb$ and on the operad $\Lie$. We denote by $\Cx_{\infty}$ the operad encoding homotopy $\Cx$-algebras. The main result of this article is the following
\begin{thmintro}[Theorem \ref{equiv1}]
There is an equivalence of categories
$$\begin{array}{ccc}
\left\{ \begin{gathered} \Cx_{\infty}\textrm{-algebras}\\ \textrm{with } \infty \textrm{-morphisms} \end{gathered} \right\} & \xrightarrow{\cong} & \left\{ \begin{gathered} \textrm{Complex structures on formal manifolds}\\ \textrm{with holomorphic maps} \end{gathered} \right\}.
\end{array}$$
\end{thmintro}
This theorem provides the algebraic essence of the geometric notion of complex manifolds, and it is given by a surprisingly simple operad. The notion of $\infty$-morphism for $\Cx_{\infty}$-algebras is analog to the notion of $\mathrm{L}_{\infty}$-morphisms (resp. $\mathrm{A}_{\infty}$-morphisms) for $\mathrm{L}_{\infty}$-algebras (resp. $\mathrm{A}_{\infty}$-algebras).\\

We are finally interested in a global version of Theorem \ref{equiv1}. Formal geometry \cite{GelfandKazhdan} gives a convenient language to provide a global description of objects defined locally in terms of coordinates. Let $M$ be a smooth manifold. The idea is to work with the space of all local coordinates systems $\Mco$. Let $x$ be a point in $M$ and $\varphi$ be a local coordinates system around $x$. Theorem \ref{equiv1} describes as a $\Cx_{\infty}$-algebra the Taylor series of a complex structure $J$ at the point $x$ in the chart given by $\varphi$. Therefore, the description of the map $J$ on $M$ is given by a collection of $\Cx_{\infty}$-algebras indexed by the points in $\Mco$. However, a collection of $\Cx_{\infty}$-algebras indexed by $\Mco$ does not necessarily correspond to a smooth endomorphism $J$. Based on the operad $\Cx_{\infty}$, we define two fiber bundles $E_{cx}(M)$ and $F_{cx}(M,\, N)$ over $M$, both of them endowed with a connection,
 such that the following theorem holds:
\begin{thmintro}[Theorem \ref{equiv2}]
Let $M$ and $N$ be two smooth manifolds. There is an equivalence of categories
$$\begin{array}{ccc}
\left\{ \begin{gathered} \textrm{Flat sections of } E_{cx}(M) \textrm{ with} \\ \textrm{flat sections of } F_{cx}(M,\, N) \end{gathered} \right\} & \xrightarrow{\cong} & \left\{ \begin{gathered} \textrm{Complex structures on } M \textrm{ with}\\ \textrm{holomorphic maps between } M \textrm{ and } N \end{gathered} \right\}.
\end{array}$$
\end{thmintro}\

\subsection*{Layout}

Since this article brings together differential geometry and operad theory, we recall quickly definitions and results from the two domains in Sections \ref{section1} and \ref{operadicinterpretation}. More precisely, we recall definitions and notations related to complex geometry in Section \ref{section1} and we fix notations for operads in Section \ref{operadicinterpretation}. In this second section, we also introduce the operad $\Cx$ and we prove that it is a Koszul operad. Moreover, we describe the algebras over the operad $\Cx_{\infty}$ and the associated $\infty$-morphisms. In Section \ref{section3}, we prove Theorem \ref{equiv1} relating $\Cx_{\infty}$-algebras and formal complex structures. The smooth version of this theorem is Theorem \ref{equiv2} and it is detailed and proved in Section \ref{section4}. There are two appendices: the first one provides an explicit decomposition map for the Koszul dual cooperad associated to Lie algebras in degree 1, and the second one explains the theory of distributive laws for cooperads.

\subsection*{Acknowledgments}

This long-term job benefits from a lot of persons. I would like to thank Sergei Merkulov who suggested this work to me and I got from him and from Henrik Strohmayer several explanations related to operad profiles. I am grateful to Mathieu Anel, Damien Calaque and Yael Frégier, as well as many of my colleagues in Toulouse: Marcello Bernardara, Denis-Charles Cisinski, Eveline Legendre, Joseph Tapia, for fruitful discussions. I wish to thank Bruno Vallette for his comments on the first version of this preprint.

I also would like to thank the Max-Planck Institut für Mathematiks, the Institut des Hautes Études Scientifiques and the Isaac Newton Institute for excellent working conditions during the stays that I have spent there.

\subsection*{Notations}

In this paper, we work over the field of real numbers $\Kb = \Rb$, except in the appendices where $\Kb$ can be any field of characteristic $0$. We use the symbol $\otimes$ for the tensor product (over $\Kb$), the symbol $\odot$ for the symmetric tensor product and the symbol $\wedge$ for the anti-symmetric product. Let $V := \{ V^n \}_{n \in \Zb}$ be a cohomologically graded vector space. We consider a one-dimensional vector space $s \Kb$ spanned by an element ``$s$" of cohomological degree $-1$ (or equivalently, of homological degree $1$ but we will only speak about cohomological degrees in this article). By definition, the \emph{cohomological suspension} of $V$ is $s^{-1}V := s^{-1}\Kb \otimes V$. It corresponds also to the shifted cochain complex $V[-1]$, so that $(s^{-1}V)^n = (V[-1])^n = V^{n-1}$. Similarly, we get the \emph{cohomological desuspension} $sV = V[1]$. The composition of a composable pair of morphisms $(f,\, g)$ is denoted by $g \cdot f$. The maps appearing in this article depend on several variables, say $f$ depends on $x$ and $t$. We denote by $d_x f$, resp. $d_t f$, the partial differential of $f$ with respect to the variables $x$, resp. $t$.

\subsection*{Conventions}

In all the paper, the manifolds and the formal manifolds are assumed to be finite dimensional. Throughout the paper, we use the Einstein summation convention, i.e. we always sum over repeated upper and lower indices. For instance, $P^{a} \partial_{a}$ means $\sum_{a} P^{a} \partial_{a}$. In the paper, we consider differential graded (dg for short) vector space that are cohomologically graded. Therefore, we assume that the differential is of degree $+1$. Moreover, we use the Koszul sign convention saying that in a commutative algebra $a_{1} a_{2} = (-1)^{|a_{1}||a_{2}|} a_{2} a_{1}$.

\section{Complex structures}\label{section1}

We remind notations and general facts concerning complex structures and formal complex structures. We also describe the holomorphic maps associated to each context.

\subsection{Complex structures}

Let $M$ be a paracompact smooth manifold. We denote by $T^*M$ its \emph{cotangent bundle} and by $TM$ its \emph{tangent bundle}. The Lie bracket on vector fields induces a symmetric product of degree $1$ on the shifted tangent bundle
$$
[-,\, -] : TM[1] \odot TM[1] \rightarrow TM[1].
$$
\begin{dei}
An \emph{almost complex structure} on $M$ is an endomorphism $J : TM[1] \rightarrow TM[1]$ satisfying $J^{2} = -\Id$.
\end{dei}
To such an endomorphism, we can associate its \emph{Nijenhuis torsion} $\N_{J} : TM[1]^{\odot 2} \rightarrow TM[1]$ given by
$$
\N_{J} (X,\, Y) := J^{2}[X,\, Y] + [JX,\, JY] - J[X,\, JY] - J[JX,\, Y].
$$
\begin{thm}[Newlander--Nirenberg \cite{NewlanderNirenberg}]
Let $J$ be an almost complex structure on $M$. The endomorphism $J$ is a complexe structure if, and only if, its Nijenhuis torsion vanishes.
\end{thm}

We denote by $\O_{M} = \C^{\infty}_{M}$ the \emph{sheaf of smooth functions} on $M$. The sheaf of $\O_M$-modules of \emph{differential $1$-forms} is given by the sheaf of sections of the shifted cotangent bundle $\Omega^{1}_{M} := \Gamma T^{*}M[-1]$ and the graded symmetric algebra $\Omega^\bullet_{M} := S^\bullet (\Omega^1_{M})$ (over the sheaf $\C^{\infty}_{M}$) has its (graded-)\-symmetric product given by the wedge product and is called the \emph{de Rham algebra}. The de Rham algebra is usually seen as an antisymmetric algebra since $\Omega^{1}_{M}$ is in degree $1$. The de Rham differential $d_{DR} : \C^{\infty}_{M} = \Omega^{0}_{M} \rightarrow \Omega^{1}_{M}$ defined by $d_{DR}f(X) = X(f)$ extends to the de Rham algebra in order to get a differential graded algebra.\\

To any \emph{vector form} $F$ in $\Omega^r_{M} \otimes_{\O_{M}} T_M[1]$, where $T_M$ is the $\O_M$-module of sections of the tangent bundle $TM$, we associate two derivations: the \emph{interior product} $i_F : \Omega^\bullet_{M} \rightarrow \Omega^{\bullet + r-1}_{M}$ defined by\\
$i_F(\omega)(Y_1 \wedge \cdots \wedge Y_{r+s-1}) :=$
$$
\frac{1}{r! (s-1)!}\sum_{\sigma \in \Sb_{r+ s-1}} \sgn \sigma \times \omega \left(F(Y_{\sigma(1)} \wedge \cdots \wedge Y_{\sigma(r)}) \wedge Y_{\sigma(r+1)} \wedge \cdots \wedge Y_{\sigma(r+s-1)}\right),
$$
for $\omega \in \Omega^s_{M}$ and $Y_1 \wedge \cdots \wedge Y_{r+s-1} \in S^{r+s-1} (T_M[1])$ and $\Sb_{r+s-1}$ is the group of permutations on $r+s-1$ elements, and the \emph{Nijenhuis--Lie derivative} $d_{F} : \Omega^{\bullet}_{M} \rightarrow \Omega^{\bullet + r}_{M}$ defined by
$$
d_{F} := i_{F} \cdot d - (-1)^{r-1} d \cdot i_{F}.
$$
The \emph{Fr\"olicher--Nijenhuis bracket} is a Lie bracket
$$
[- ,\, - ]_{\textrm{F-N}} : (\Omega^r_{M} \otimes_{\O_{M}} T_M[1]) \otimes (\Omega^s_{M} \otimes_{\O_{M}} T_M[1]) \rightarrow \Omega^{r+s}_{M} \otimes_{\O_{M}} T_M[1]
$$
given, for  $F = \varphi \otimes X \in \Omega^r_{M} \otimes_{\O_{M}} T_M[1]$ and $G = \psi \otimes Y \in \Omega^s_{M} \otimes_{\O_{M}} T_M[1]$, explicitly by\\
$[F,\, G]_{\textrm{F-N}} :=$
$$
\varphi \wedge \psi \otimes [X,\, Y] + \varphi \wedge d_{X}(\psi) \otimes Y - d_{Y}(\varphi) \wedge \psi \otimes X + (-1)^{r} \left( d\varphi \wedge i_{X}(\psi) \otimes Y + i_{Y}(\varphi) \wedge d\psi \otimes Y\right).
$$\

An endomorphism $J : T_M[1] \rightarrow T_M[1]$ can be seen as an element $J$ in $\Omega^{1}_{M} \otimes_{\O_{M}} T_M[1]$. Similarly, we get $[J,\, J]_{\textrm{F-N}} \in \Omega^{2}_{M} \otimes_{\O_{M}} T_M[1]$.

\begin{thm}[Nijenhuis \cite{Nijenhuis}]
Under the identification between elements and applications, 
$$
\N_{J} = [J,\, J]_{\textrm{F-N}}.
$$
\end{thm}
For an element $J \in \Omega^{1}_{M} \otimes_{\O_{M}} T_M[1]$, we call it \emph{integrable} if $[J,\, J]_{\textrm{F-N}} =0$. This means that a complex structure is an almost complex structure $J$ satisfying the integrability condition.\\

\begin{dei}
Let $(M,\, J)$ and $(N,\, J')$ be two complex manifolds. A smooth map $F : (M,\, J) \rightarrow (N,\, J')$ is called \emph{$J$-$J'$-holomorphic}, or \emph{holomorphic} when there is no ambiguity, if it satisfies
$$
dF \cdot J = J'(F) \cdot dF.
$$
\end{dei}

\subsection{Complex formal manifolds}

In order to work with formal power series instead of smooth functions, we will now consider \emph{formal pointed manifolds}, that is, locally ringed spaces of the form $\V_{\mathrm{for}} = \left(\{\textrm{point} \},\, \O_{\V_{\mathrm{for}}}\right)$, where the sheaf of functions is given by power series $\O_{\V_{\mathrm{for}}} := S(V)^{*}\ (= \mathrm{Hom}(S(V),\, \Rb))$. It is the formal scheme associated to $\left(\mathrm{Spec}(S(V)),\, S(V)\right)$ and it can be interpreted as a formal neighborhood of a base point, say $0$, in the vector space $V$. Its \emph{cotangent sheaf} is the formal completion of the cotangent sheaf of $\left(\mathrm{Spec}(S(V)),\, S(V)\right)$ and is given by $\Omega_{\V_{\mathrm{for}}} \cong \O_{\V_{\mathrm{for}}} \otimes (V[1])^*$. Similarly, its \emph{tangent sheaf} is the $\O_{\V_{\mathrm{for}}}$-module $T_{\V_{\mathrm{for}}} \cong  \O_{\V_{\mathrm{for}}} \otimes V$.\\

Let $\{ e_{a} \}$ be a basis of $V$, or $\left\{ e_{a}^{V} \right\}$ when there is an ambiguity, and $\{ t^{a} \}$ be its associated dual basis, so that $\O_{\V_{\mathrm{for}}} \cong \Rb \llbracket t^{a} \rrbracket$. The shifted tangent sheaf $T_{\V_{\mathrm{for}}}[1]$ is generated as an $\O_{\V_{\mathrm{for}}}$-module by a basis $\{ \partial_{a} \}$, where we write $\partial_{a}$ for $s\frac{\hspace{-.2cm}\partial}{\partial t^{a}}$ (we remind that ``$s$'' stands for the cohomological desuspension, that is to say, $\partial_{a}$ is of degree $-1$). Finally, we denote by $\{ \gamma^{a} \}$ the dual basis of $\{ \partial_{a}\}$, that is $\gamma^a = s^{-1}dt^{a}$, where $dt^{a}$ is dual to $\frac{\hspace{-.2cm}\partial}{\partial t^{a}}$. The element $\gamma^{a}$ is therefore of (cohomological) degree $1$. We obtain that the $\O_{\V_{\mathrm{for}}}$-module of differential forms is $\Omega^{\bullet}_{\V_{\mathrm{for}}} \cong \Rb \llbracket t^{a},\, \gamma^{b} \rrbracket$.\\

A \emph{vector valued differential form}, or \emph{vector form} for short, is an element $F \in\Omega^{\bullet}_{\V_{\mathrm{for}}} \otimes_{\O_{\V_{\mathrm{for}}}} T_{\V_{\mathrm{for}}}$. It has the following form
$$
F = \sum_{p \geq 0} F_{a_{1} \cdots a_{p}}^{b}(t) \gamma^{a_{1}} \cdots \gamma^{a_{p}} \partial_{b},
$$
where $F_{a_{1} \cdots a_{p}}^{b}(t)$ is an element in $\O_{\V_{\mathrm{for}}}$, that is a power series in the $t^{a}$'s. For example, an endomorphism $J : T_{\V_{\mathrm{for}}} \rightarrow T_{\V_{\mathrm{for}}}$, given by $J(P^{a}(t) \partial_{a}) = P^{a}(t) J_{a}^{b}(t) \partial_{b}$, where $P^{a}(t)$, $J_{a}^{b}(t) \in \O_{\V_{\mathrm{for}}}$, is seen as a vector form $J = J_{a}^{b}(t) \gamma^{a} \partial_{b} \in \Omega^{1}_{\V_{\mathrm{for}}} \otimes_{\O_{\V_{\mathrm{for}}}} \T_{\V_{\mathrm{for}}}[1]$. For two vector forms $F = \sum_{p} F_{a_{1} \cdots a_{p}}^{b}(t) \gamma^{a_{1}} \cdots \gamma^{a_{p}} \partial_{b}$ and $G = \sum_{q} G_{c_{1} \cdots c_{q}}^{d}(t) \gamma^{c_{1}} \cdots \gamma^{c_{q}} \partial_{d}$, a formula for the Fr\"olicher--Nijenhuis bracket is given by\\
$[F,\, G]_{\textrm{F-N}} = \left( F_{a_{1} \cdots a_{p}}^{b}(t) \partial_{b} G_{a_{p+1} \cdots a_{p+q}}^{d}(t) - G_{a_{p+1} \cdots a_{p+q}}^{b}(t) \partial_{b} F_{a_{1} \cdots a_{p}}^{d}(t)\right.$
$$
\left. - p F_{b a_{2} \cdots a_{p}}^{d}(t) \partial_{a_{1}} G_{a_{p+1} \cdots a_{p+q}}^{b}(t) + q G_{b a_{p+2} \cdots a_{p+q}}^{d}(t) \partial_{a_{p+1}} F_{a_{1} \cdots a_{p}}^{b}(t) \right) \gamma^{a_{1}} \cdots \gamma^{a_{p+q}} \partial_{d}.
$$
For instance, for $J = J_{a}^{b}(t) \gamma^{a} \partial_{b}$, we get the formula:
\begin{equation}\label{bracket}
[J,\, J]_{\textrm{F-N}} = \left( J_{a_{1}}^{b} \partial_{b} J_{a_{2}}^{d} - J_{a_{2}}^{b} \partial_{b} J_{a_{1}}^{d} - J_{b}^{d} \partial_{a_{1}} J_{a_{2}}^{b} + J_{b}^{d} \partial_{a_{2}} J_{a_{1}}^{b} \right) \gamma^{a_{1}} \gamma^{a_{2}} \partial_{d},
\end{equation}
where we have removed the variable $t$.\\

A map of formal manifold $F : \V_{\mathrm{for}} = \left(\{\textrm{pt} \},\, \O_{\V_{\mathrm{for}}}\right) \rightarrow \W_{\mathrm{for}} = \left(\{\textrm{pt} \},\, \O_{\W_{\mathrm{for}}}\right)$ is of the form $F = F^{b}(t) e_{b}^{W}$, where $\{ e_{a}^{W} \}$ is a basis of $W$. Therefore, its differential $dF : T_{\V_{\mathrm{for}}} \rightarrow F^*T_{\W_{\mathrm{for}}}$ is given by
$$
dF = \partial_{a} F^{b}(t) \gamma^{a} \partial_{b}.
$$
It follows that a holomorphic map between two complex formal manifolds is a map $F : (\V_{\mathrm{for}},\, J) \rightarrow (\W_{\mathrm{for}},\, J')$ satisfying
$$
\left( \partial_{c} F^{d}(t) \gamma^{c} \partial_{d} \right) \cdot \left( J_{a}^{b}(t) \gamma^{a} \partial_{b} \right) = \left( {J'}_{c}^{d}(F(t)) \gamma^{c} \partial_{d} \right) \cdot \left( \partial_{a} F^{b}(t) \gamma^{a} \partial_{b} \right),
$$
that is to say,
\begin{equation}\label{holom}
\partial_{b} F^{d}(t) J_{a}^{b}(t) \gamma^{a} \partial_{d} = {J'}_{b}^{d}(F(t))  \partial_{a} F^{b}(t) \gamma^{a} \partial_{d}.
\end{equation}

\section{Operadic interpretation}\label{operadicinterpretation}

In this section, we define the operad $\Cx$ and we show that it is a Koszul operad. In order to prove this fact, we apply the distributive laws theory recalled in Appendix \ref{distri} and we make the curved Koszul dual cooperad $\Cx^{\ash}$ explicit. We finally described the homotopy algebras associated to the operad $\Cx$ and the corresponding homotopy morphisms.

\subsection{Definitions and notations}

We refer to the book written by Loday and Vallette \cite{LodayVallette} for definitions about operads. However, we consider cohomological grading and therefore differentials have degree $+1$ and the homological suspension will be the cohomological desuspension. Let $\P = \T(E)/(R)$ be a quadratic operad, where $\T(E)$ is the free operad generated by the $\Sb$-module $E$ and $(R)$ is the ideal generated by the $\Sb$-module $R$. Its Koszul dual cooperad $\P^{\ash}$ is a subcooperad of the cofree cooperad $\T^{c}(sE)$ on $sE$, and is defined by the universal property dual to the universal property defining the quotient, for the cogenerators $sE$ and the coideal cogenerated by $s^{2}R$. We use the notation $\P^{\ash} = \C(sE\, ;\, s^2 R)$.  We denote by $\Im$ the $\Sb$-module $(0,\, \Rb ,\, 0,\, \ldots)$ and by $M \circ_{(1)} N$ the infinitesimal composite of two $\Sb$-modules $M$ and $N$. For a coaugmented cooperad $(\C,\, \Delta_{\C})$, we denote by $\Delta_{\C}^{(1)}$ the \emph{infinitesimal} decomposition map $\C \rightarrow \C \circ_{(1)} \C$ (see Section 6.1 in \cite{LodayVallette} for precise definitions) and by $\bar{\Delta}_{\C}^{(1)}$ the \emph{reduced} infinitesimal decomposition map $\C \xrightarrow{\Delta_{\C}^{(1)}} \C \circ_{(1)} \C \twoheadrightarrow \overline{\C} \circ_{(1)} \overline{\C}$.\\

For instance, we denote by $\Lie_{1}$ the operad encoding Lie algebras with a bracket of cohomological degree $1$, that is to say,
$$
\Lie_{1} := \T \left(s^{-1} E_{L} \right)/\left(s^{-2} R_{L} \right),
$$
where $E_{L}$ is the $\Sb$-module generated in arity $2$ by a symmetric element $\lie$ and $R_{L}$ is the $\Sb$-module generated in arity $3$ by the Jacobi relation: $\jac$. In this case, we get $\Lie_{1}^{\ash} = \C (E_{L}\, ;\, R_{L}) = \Com^{*}$. Here, we denote by $\Com$ the operad encoding commutative algebra (with a product of degree $0$) and therefore, the cooperad obtained by dualizing arity-wise, $\Com^{*}$, is the cooperad encoding cocommutative coalgebras. The Koszul dual cooperad $\Lie_{1}^{\ash}$ is 1-dimensional in each arity, $\Lie_{1}^{\ash}(n) \cong \Rb \cdot \overline{l}_{n}^{c}$, where $\overline{l}_{n}^{c}$ is an element of degree $0$ on which $\Sb_{n}$ acts trivially, and the infinitesimal decomposition map on $\Lie_{1}^{\ash}$ is given by
$$
\Delta^{(1)}_{\Lie_{1}^{\ash}}(\overline{l}_{n}^{c}) = \sum_{\substack{p+q=n+1\\ p,\, q \geq 1}} \sum_{\sigma \in Sh_{q,\, p-1}^{-1}} (\overline{l}_{p}^{c} \circ_{1} \overline{l}_{q}^{c})^{\sigma},
$$
where $Sh_{q,\, p-1}^{-1}$ is the set of $(q,\, p-1)$-unshuffles, that is, inverses of $(q,\, p-1)$-shuffles. We refer to Appendix \ref{Lie} for more details.

\subsection{Toward the operad profile of complex structures}

In Section \ref{section1}, we have seen that an element $J \in \Omega^{1}_{\V} \otimes_{\O_{\V}} T_{\V}[1]$ is a complex structure on a formal pointed manifold $\V_{\mathrm{for}} = \left(\{\textrm{pt} \},\, \O_{\V_{\mathrm{for}}}\right)$ if and only if the following two equations are satisfied:
\begin{equation}\label{rel} 
\left\{ \begin{split}
J^{2} +\Id = 0,\\
{[J,\, J]_{\textrm{F-N}}} = 0.
\end{split} \right.
\end{equation}
The element $J$ has the form $J = J_{a}^{b}(t) \gamma^{a} \partial_{b}$, where $J_{a}^{b}(t) = \sum J_{c_{1},\, \ldots ,\, c_{n};\, a}^{b} t^{c_{1}} \cdots t^{c_{n}}$. The two previous equations give the relations that the coefficients $J_{c_{1},\, \ldots ,\, c_{n};\, a}^{b}$ have to satisfied. Our aim is now to describe such a structure as an algebraic data.\\

We remark that every map $f : V^{\otimes n+1} \rightarrow V$ is given by its values on a basis
$$
f(e_{c_{1}} \otimes \cdots \otimes e_{c_{n}} \otimes e_{a}) = f_{c_{1},\, \ldots ,\, c_{n},\, a}^{b} e_{b}
$$
and we can therefore think about $J$ as a sum of applications $V^{\otimes n+1} \rightarrow V$ satisfying relations. In order to make the problem simpler, we begin with a complex structure $J$ equal to its constant part $J_{a}^{b}(0) = J_{a}^{b}$, that is to say, since a basis of $V$ is fixed, $J$ is an endomorphism $V \rightarrow V$. In that case, the Nijenhuis torsion is trivially equal to $0$ since $J$ is constant. Therefore we can see such a simple complex structure as a representation of the unital associative algebra of complex numbers $\Cb$ in $V$, that is, an algebra morphism
$$
\Cb \rightarrow \End(V) := \Hom (V,\, V).
$$
In operadic terms, we say that the vector space $V$ is an algebra over the (nonsymmetric) operad
$$
\underline{\Cb} := (0,\, \Cb,\, 0,\, \ldots) \cong \left. \T(\jnsi) \middle/ \left( \jnrsi \right) \right. .
$$
To pass from algebra to geometry, we have to handle the question of a change of the coordinates system. Such a data for a formal pointed manifold $\V_{\mathrm{for}}$ is given by a \emph{formal diffeomorphism} $\phi$ of $\V_{\mathrm{for}}$, that is, a power series with values in $V$ fixing the base point and which linear part is invertible. In other words, it is an application of vector spaces $\phi : \overline{S}(V) \cong \Lie_{1}^{\ash}(V) \rightarrow V$, such that the linear part $\phi_{1} : V \rightarrow V$ is an isomorphism. It can be seen equivalently as an $\infty$-isomorphism of the trivial $\Lie_{1}$-algebra $V$ to itself. (Here $V$ is considered concentrated in cohomological degree $0$.) We therefore add to our algebraic data $\underline{\Cb}$ a generator $\slie$ and the suspended Jacobi relation which encode Lie-brackets in cohomological degree $1$. In this way, formal diffeomorphisms will appear in the homotopy theory of our algebraic data. The condition for a change of coordinates can be written $d\phi \cdot J = J'(\phi) \cdot d\phi$. We focus here our attention on the generators of the algebraic data $\jnsi$ and $\slie$. We consider two representations $J_a^b \gamma^a \partial_b$ and ${J'}_a^b \gamma^a \partial_b$ of $\jnsi$ and a representation of $\slie$ corresponding to the quadratic part $q_{a b}^c t^a t^b e_c$ of $\phi$. The condition $d\phi \cdot J = J'(\phi) \cdot d\phi$ implies the algebraic equations
$$q_{b d}^c J_a^d = {J'}_d^c q_{b a}^d.$$
As we will see in Proposition \ref{infinitymorphism} and Corollary \ref{corinfinitymorph}, these relations correspond to the representation for morphisms of the operadic relation
$$
\sjlr = \sljr,
$$
and it will therefore appear in the algebraic data.

\begin{rem}
This reasoning is based on the fact that the category of formal manifolds with formal smooth maps is equivalent to the category of $(\Lie_{1})_{\infty}$-algebras with $\infty$-morphisms.
\end{rem}

\begin{dei}
We define the operad
$$
\Cx := \left. \T\left( \jn ,\, \slie \right) \middle/ \left( \qjnr ,\, \sjlr - \sljr ,\, \scriptscriptstyle{s^{-2}}\left( \jac \right) \right) \right. .
$$
\end{dei}

\begin{rems} \leavevmode
\begin{itemize}
\item In this notation, we assume that the element $\sjlrt - \sljrt$ is a relation in $\Cx$.
\item In the description of $\underline{\Cb}$ and in the sequel, since there is no ambiguity, we sometimes omit the 1 on elements in arity 1 as $\jn$ and write for instance $\jnsi$.
\end{itemize}
\end{rems}

We prove that $\Cx$ is the operad profile of complex structures in the rest of Section \ref{operadicinterpretation} and in Section \ref{section3}.

\subsection{Distributive laws and the Koszul dual cooperad}

The operad $\Cx$ is not a homogeneous quadratic operad because of the fact that the relation $\jnrsi$ involves quadratic and constant terms. We use the theory developed by Hirsh and the author in \cite{HirshMilles} to find an explicit cofibrant replacement of $\Cx$. The quadratic operad $\qCx$ associated to $\Cx$ has the following presentation
$$
\qCx := \left. \T\left( \jn ,\, \slie \right) \middle/ \left( \qjnr ,\, \sjlr - \sljr ,\, \scriptscriptstyle{s^{-2}}\left( \jac \right) \right) \right. .
$$
By means of distributive laws, we make the Koszul dual cooperad $\qCx^{\ash}$ explicit and we prove that the operad $\qCx$ is Koszul. We obtain finally that $\Cx$ is Koszul. We refer to Loday-Vallette \cite{LodayVallette}, Section 8.6 for definitions and notations on distributive laws.\\

The quadratic operad (algebra) associated to $\underline{\Cb}$ is
$$
\qCb := \T(\jnsi)/\left( \qjnrsi \right) = \Rb \left[ \jnsi \right].
$$
The Koszul dual cooperad is given by $\displaystyle{\qCb^{\ash} = \Rb \left[\sojnsi \right] = \bigoplus_{n \geq 0} \Rb\, \overline{\imath}_{n}^{c}}$, where we denote $s\jnsi$ by $\sojnsi$ and where $\overline{\imath}_{n}^{c}$ is an element of degree $-n$ and has the following (infinitesimal) decomposition map
$$
\Delta_{\qCb^{\ash}}(\overline{\imath}_{n}^{c}) = \sum_{k =0}^{n} (\overline{\imath}_{k}^{c}\, ;\, \overline{\imath}_{n-k}^{c}) \in \qCb^{\ash} \circ \qCb^{\ash} \cong \qCb^{\ash} \otimes \qCb^{\ash}.
$$
Between the operads $\Lie_{1}$ and $\qCb$, we define the rewriting rule
$$
\lambda : (s^{-1}E_{L}) \circ_{(1)} \left( \Rb \jnsi \right) \rightarrow \left( \Rb \jnsi \right) \circ_{(1)} (s^{-1}E_{L})
$$
by the $\Sb_{2}$-equivariant map sending $\sjlr$ to $\sljr$. We refer to Appendix~\ref{distri} for general facts on distributive laws on cooperads. We use the notations $\qCb \vee_{\lambda} \Lie_{1} = \qCx$ and
\begin{equation*}
\begin{split}
\qCb^{\ash} \vee^{\lambda} \Lie_{1}^{\ash} & = \qCx^{\ash}\\
& = \C \left(s\left( \jn ,\, \slie \right)\, ;\, s^{2}\left( \qjnr ,\, \sjlr - \sljr ,\, \scriptscriptstyle{s^{-2}}\left( \jac \right) \right) \right)\\
& = \C \left( \sojn ,\, \lie\, ;\, \soqjnr ,\, \sojlr + \soljr ,\, \jac \right).
\end{split}
\end{equation*}
There is the change of sign $s^2 \left(\sjlr - \sljr \right) = -\ \sojlr - \soljr$ because of the Koszul sign rule.

\begin{lem}\label{lambda}
The injection \emph{$i_{1} : \qCb^{\ash} \vee^{\lambda} \Lie_{1}^{\ash} = \qCx^{\ash} \hookrightarrow \Lie_{1}^{\ash} \circ \qCb^{\ash}$} (defined in Appendix \ref{naturali}) is an isomorphism. Therefore, the morphism of $\Sb$-modules $\lambda$ induces a distributive law of operads $\Lambda : \Lie_{1} \circ \qCb \rightarrow \qCb \circ \Lie_{1}$ and a distributive law of cooperads \emph{$\Lambda^{c} : \Lie_{1}^{\ash} \circ \qCb^{\ash} \rightarrow \qCb^{\ash} \circ \Lie_{1}^{\ash}$}. The distributive law of cooperads $\Lambda^{c}$ is given explicitly by
$$
\Lambda^{c}\left((\overline{l}_{n}^{c}\, ;\, \overline{\imath}_{k_{1}}^{c},\, \ldots ,\, \overline{\imath}_{k_{n}}^{c})\right) = \left(\sum_{\sigma \in Sh_{k_{1},\, \ldots ,\, k_{n}}} \sgn \sigma \right) (\overline{\imath}_{k_{1} + \cdots + k_{n}}^{c}\, ;\, \overline{l}_{n}^{c}),
$$
where $Sh_{k_{1},\, \ldots ,\, k_{n}}$ is the set of $(k_{1},\, \ldots ,\, k_{n})$-shuffles.
\end{lem}

\begin{rem}
Pictorially, we get the following formula:
$$
\Lambda^{c} \left( {\tiny \vcenter{\xymatrix@M=0pt@R=8pt@C=10pt{
    & & & & &\\
    & \circ \ar@{-}[u] \ar@{}[]^{k_{2}} & \cdots & & \circ \ar@{}[]^{k_{n-1}} \ar@{-}[u] & \circ \ar@{}[]_{k_{n}} \ar@{-}[u] \\
    \ar@{-}[drrr] & \ar@{-}[u] \ar@{-}[drr] & \ar@{-}[dr] & \ar@{-}[d] & \ar@{-}[dl] \ar@{-}[u] & \ar@{-}[dll] \ar@{-}[u] \\
    & & &  \ar@{-}[d] & &\\
    & & & & &
    }}} \right) = \alpha_{k_1,\, \ldots ,\, k_n} {\tiny \vcenter{\xymatrix@M=0pt@R=8pt@C=10pt{
    & & & & &\\
    \ar@{-}[drrr] & \ar@{-}[drr] & \ar@{-}[dr] & \ar@{-}[d] & \ar@{-}[dl] & \ar@{-}[dll] \\
    & & &  \ar@{-}[d] & &\\
    &&& \circ \ar@{}[]^{k_{1}+\cdots +k_{n}} \ar@{-}[d] &&\\
    & & & & &
    }}},
$$
where $\alpha_{k_1,\, \ldots ,\, k_n}$ is a coefficient in $\Zb$ and where $\scriptstyle{k}\! \! \sojnsi$ stands for $\bar{\imath}_k^c$.
\end{rem}

\begin{pf}
Theorem \ref{diamond} says that it is enough to prove that $p$ is injective in weight $-3$, where the weight is given by the opposite of the number of generators. The last remark of Section \ref{remdiamond} ensures that we can equivalently prove the surjectivity of $i_{1}$ in weight $-3$. In order to make the reader more familiar with $\qCx^{\ash}$, we will make the elements of weight $-3$ in $\qCx^{\ash}$ explicit. The map $i_{1}$ is trivially surjective onto $\Lie_{1}^{\ash} \circ \Im$ and onto $\Im \circ \qCb^{\ash}$. The other terms of weight $-3$ in $\qCx^{\ash} = \qCb^{\ash} \vee^{\lambda} \Lie_{1}^{\ash}$ are
$$
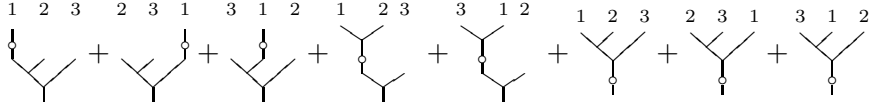
\jjac
$$
and the $\Sb_{3}$-module generated by it,
$$
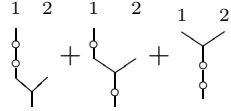
\jjlr
$$
and the $\Sb_{2}$-module generated by it, and finally the element
$$
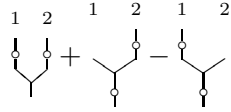
\jljr
$$
on which $\Sb_{2}$ acts by signature. This proves that $i_{1}$ is surjective in weight $-3$.

To describe $\Lambda^{c}$, we first remark that the formula is true for $n=1$, for $k_{1}+\cdots +k_{n} = 0$ and, for $n=2$ and $k_{1}+k_{2}=1$. Then, we make use of diagrams (I) and (II) in Appendix \ref{distrilaw} to prove the general formula. Let $k \in \Nb^{*}$. We assume that the formula is true for any $m \leq n$ and $k_{1}+\cdots +k_{m} \leq k-1$ and we use diagram (II) to prove it for $n$ and $k_{1}+\cdots +k_{n} = k$. We have
\begin{multline*}
(\Id \circ \Lambda^{c}) \cdot (\Lambda^{c} \circ \Id) \cdot (\Id \circ \Delta_{\qCb^{\ash}})(\overline{l}_{n}^{c}\, ;\, \overline{\imath}_{k_{1}}^{c},\, \ldots ,\, \overline{\imath}_{k_{n}}^{c})\\
\shoveleft{= (\Id \circ \Lambda^{c}) \cdot (\Lambda^{c} \circ \Id) \left(\sum_{k_{j}' + k_{j}'' = k_{j},\, k'_{j},\, k''_{j} \geq 0} \varepsilon_{k_{j}',\, k''_{j}} (\overline{l}_{n}^{c}\, ;\, \overline{\imath}_{k'_{1}}^{c},\, \ldots ,\, \overline{\imath}_{k'_{n}}^{c}\, ;\, \overline{\imath}_{k''_{1}}^{c},\, \ldots ,\, \overline{\imath}_{k''_{n}}^{c}) \right)}\\
\shoveleft{= (\Id \circ \Lambda^{c}) \Bigg((T\, ;\, \Id,\, \ldots ,\, \Id) + (\Id\, ;\, \overline{l}_{n}^{c}\, ;\, \overline{\imath}_{k_{1}}^{c},\, \ldots ,\, \overline{\imath}_{k_{n}}^{c})}\\
\shoveright{\left. + \sum_{\substack{k_{j}' + k_{j}'' = k_{j}\\ 1 \leq k'_{1}+\cdots +k'_{n} \leq n-1}} \varepsilon_{k_{j}',\, k''_{j}} \left( \sum_{\sigma' \in Sh_{k'_{1},\, \cdots ,\, k'_{n}}} \sgn \sigma' \right) (\overline{\imath}_{k'_{1}+ \cdots +k'_{n}}^{c}\, ;\, \overline{l}_{n}^{c}\, ;\, \overline{\imath}_{k''_{1}}^{c},\, \ldots ,\, \overline{\imath}_{k''_{n}}^{c}) \right)}\\
\shoveleft{= (T'\, ;\, \Id,\, \ldots ,\, \Id\, ;\, T'') + (\Id\, ;\, T' \, ;\, T'')}\\
+ \sum_{\substack{k_{j}' + k_{j}'' = k_{j}\\ 1 \leq k'_{1}+\cdots +k'_{n} \leq n-1}} \varepsilon_{k_{j}',\, k''_{j}} \sum_{\sigma' \in Sh_{k'_{1},\, \cdots ,\,k'_{n}}} \sgn \sigma' \sum_{\sigma'' \in Sh_{k''_{1},\, \cdots ,\, k''_{n}}} \sgn \sigma'' \cdot (\overline{\imath}_{k'_{1}+ \cdots +k'_{n}}^{c}\, ;\, \overline{\imath}_{k''_{1}+ \cdots +k''_{n}}^{c}\, ;\, \overline{l}_{n}^{c})
\end{multline*}
where $\varepsilon_{k'_{j},\, k''_{j}} = (-1)^{\sum_{i=1}^{n-1} k_{i}''(k_{i+1}' + \cdots + k_{n}')}$ and
$$
T = (T'\, ;\, T'') = \Lambda^{c}(\overline{l}_{n}^{c}\, ;\, \overline{\imath}_{k_{1}}^{c},\, \ldots ,\, \overline{\imath}_{k_{n}}^{c}) =: \alpha_{k_{1},\, \ldots ,\, k_{n}} (\overline{\imath}_{k_{1}+\cdots +k_{n}}^{c} \, ;\, \overline{l}_{n}^{c}).
$$
On the other side
$$
(\Delta_{\qCb^{\ash}} \circ \Id) \cdot \Lambda^{c} (\overline{l}_{n}^{c}\, ;\, \overline{\imath}_{k_{1}}^{c},\, \ldots ,\, \overline{\imath}_{k_{n}}^{c}) = (\Delta_{\qCb^{\ash}} \circ \Id) (T) = \alpha_{k_{1},\, \ldots ,\, k_{n}} \sum_{k' + k'' = k} (\overline{\imath}_{k'}^{c}\, ;\, \overline{\imath}_{k''}^{c}\, ;\, \overline{l}_{n}^{c}).
$$
To prove that $\alpha_{k_{1},\, \ldots ,\, k_{n}} = \sum_{\sigma \in Sh_{k_{1},\, \ldots ,\, k_{n}}} \sgn \sigma$, it is then enough to check that for $1 \leq k \leq n-1$, we have
$$
\sum_{\substack{k_{j}' + k_{j}'' = k_{j}\\ k'_{1}+\cdots +k'_{n} = k}} \varepsilon_{k_{j}',\, k''_{j}} \left( \sum_{\sigma' \in Sh_{k'_{1},\, \cdots ,\,k'_{n}}} \sgn \sigma' \right) \left(\sum_{\sigma'' \in Sh_{k''_{1},\, \cdots ,\, k''_{n}}} \sgn \sigma'' \right) = \sum_{\sigma \in Sh_{k_{1},\, \ldots ,\, k_{n}}} \sgn \sigma.
$$
This follows from the fact that, when $k$ is fixed, there is a unique decomposition of a $(k_{1},\, \ldots ,\, k_{n})$-shuffle $\sigma$ in the following manner: some positive $k'_{j} \leq k_{j}$ such that $k'_{1}+\cdots +k'_{n} = k$, a $(k'_{1},\, \ldots ,\, k'_{n})$-shuffle $\sigma'$, a $(k''_{1},\, \ldots ,\, k''_{n})$-shuffle $\sigma''$ and a permutation $\tau \in \Sb_{k_{1}+\cdots +k_{n}}$ that sends
$$
(1,\, \ldots ,\, k'_{1},\, \ldots ,\, k_{1},\, \ldots ,\, k_{1}+\cdots + k_{n-1},\, \ldots ,\,  k_{1}+\cdots + k_{n-1} + k'_{n},\, \ldots ,\, k_{1}+\cdots + k_{n})
$$
to
$$
(1,\, \ldots ,\, k'_{1},\, k_{1}+1,\, \ldots ,\, k_{1}+k'_{2},\, \ldots ,\,  k_{1}+\cdots + k_{n-1} + k'_{n},\, k'_{1}+1,\, \ldots ,\, k_{1},\, \ldots ,\, k_{1}+\cdots + k_{n}).
$$
The signature of $\tau$ is precisely $\varepsilon_{k'_{j},\, k''_{j}}$.

Similarly, diagram (I) proves that the formula is true for $n$ and $k_{1}+\cdots +k_{n} = k$ whenever it is true for $m \leq n-1$ and $k_{1}+\cdots +k_{m} \leq k$. This concludes the proof.
$\cqfd$
\end{pf}

\subsection{Koszulity}

In this section, following the theory developed in \cite{HirshMilles} ,we compute the Koszul dual curved cooperad associated to the operad $\Cx$ and we prove that $\Cx$ is Koszul.

\begin{prop}\label{Koszul}
The operad $\qCx$ is Koszul and its Koszul dual cooperad \emph{$\qCx^{\ash}$} is given by \emph{
$$
\qCx^{\ash} \cong \Lie_{1}^{\ash} \circ \qCb^{\ash} = \Lie_{1}^{\ash} \circ \Rb \left[\sojnsi \right],
$$}
with the infinitesimal decomposition map given on generators $\overline{\jmath}^c_{k_1,\, \ldots ,\, k_n} := (\overline{l}_n^c\, ;\, \overline{\imath}_{k_1}^c,\, \ldots ,\, \overline{\imath}^c_{k_n})$ by
\begin{multline*}
\Delta_{\Lambda^{c}}^{(1)} \left( \overline{\jmath}^c_{k_1,\, \ldots ,\, k_n}\right) =\\
\sum_{\substack{p+q=n+1\\ p,\, q \geq 1}} \sum_{\sigma \in Sh_{q,\, p-1}} \sum_{\substack{k'_{\sigma(j)}+k''_{\sigma(j)} = k_{\sigma(j)}\\ k'_{\sigma(j)} = k_{\sigma(j)} \textrm{ for } j > q\\ k' = k'_{\sigma(1)}+ \cdots +k'_{\sigma(q)}}} \alpha^{\sigma}_{k'_{j},\, k''_{j}} \times \left( \overline{\jmath}_{k',\, k_{\sigma(q+1)},\, \ldots ,\, k_{\sigma(n)}}^{c} \circ_{1} \overline{\jmath}_{k''_{\sigma(1)},\, \ldots ,\, k''_{\sigma(q)}}^{c} \right)^{\sigma^{-1}},
\end{multline*}
where $\alpha^{\sigma}_{k'_{j},\, k''_{j}} := \sgn_{k_{1},\, \ldots ,\, k_{n}}\sigma \times \varepsilon_{k'_{j},\, k''_{j}}^{\sigma} \times \alpha_{k'_{\sigma(1)},\, \ldots ,\, k'_{\sigma(q)}}$ with
$$
\left\{ \begin{array}{lcl}
\varepsilon_{k'_{j},\, k''_{j}}^{\sigma} & := & (-1)^{\sum_{i=1}^{n} k_{\sigma(i)}''(k_{\sigma(i+1)}' + \cdots + k_{\sigma(n)}')}\\
\alpha_{k'_{1},\, \ldots ,\, k'_{q}} & := & \sum_{\sigma' \in Sh_{k'_{1},\, \cdots ,\,k'_{q}}} \sgn \sigma',\, \textrm{with convention } \alpha_{0,\, \ldots ,\, 0} := 1,
\end{array}\right.
$$
and $\sgn_{k_{1},\, \ldots ,\, k_{n}}\sigma$ is the signature of the restriction of $\sigma$ to the indices $j$ such that $k_{j}$ is odd (after relabeling the remaining $\sigma(j)$ in a way that the order of the $\sigma(j)$ does not change). Moreover, the full decomposition map is given by
\begin{multline*}
\Delta_{\Lambda^{c}} \left( \overline{\jmath}^c_{k_1,\, \ldots ,\, k_n}\right) =\\
\sum_{\substack{\{q_{1},\, \ldots ,\, q_{p} \}\\ q_{1}+\cdots +q_{p}= n}} \sum_{\sigma \in Sh_{q_{1},\, \ldots ,\, q_{p}}} \sum_{k'_{j}+k''_{j} = k_{j}} \beta^{\sigma}_{k'_{j},\, k''_{j}} \times \left( \overline{\jmath}_{l'_1,\, \ldots ,\, l'_p}^{c}\, ;\, \overline{\jmath}_{k''_{\sigma(1)},\, \ldots ,\, k''_{\sigma(q_1)}}^{c},\, \overline{\jmath}_{k''_{\sigma(q_1+1)},\, \ldots ,\, k''_{\sigma(q_1+q_2)}}^{c} ,\, \ldots \right)^{\sigma^{-1}},
\end{multline*}
where $l'_i := k'_{\sigma(q_1+\cdots +q_{i-1}+1)}+ \cdots +k'_{\sigma(q_1+\cdots +q_i)}$ and
$$
\beta^{\sigma}_{k'_{j},\, k''_{j}} := \sgn_{k_{1},\, \ldots ,\, k_{n}}\sigma \times \varepsilon_{k'_{j},\, k''_{j}}^{\sigma} \times \prod_{i=1}^p \alpha_{k'_{\sigma(q_1+\cdots +q_{i-1}+1)},\, \ldots ,\, k'_{\sigma(q_1+\cdots +q_i)}}.
$$
\end{prop}

\begin{rem}
The second sum really runs on shuffles and not on unshuffles since the unshuffle permutation $\sigma$ acts as a left-action given by $\sigma \cdot (\overline{\imath}_{k_{1}}^{c},\, \ldots ,\, \overline{\imath}_{k_{n}}^{c}) = \left(\overline{\imath}_{k_{\sigma^{-1}(1)}}^{c},\, \ldots ,\, \overline{\imath}_{k_{\sigma^{-1}(n)}}^{c}\right)^{\sigma}$.
\end{rem}

\begin{pf}
We have seen in the proof of Lemma \ref{lambda} that we can apply Theorem \ref{diamond}, this gives that the operad $\qCx^{\ash}$ is Koszul. It is then enough to prove the formula for the infinitesimal decomposition map by means of Proposition \ref{decomposition map}. The decomposition map on $\qCx^{\ash} \cong \Lie_{1}^{\ash} \circ \qCb^{\ash}$ is given by $\Delta_{\Lambda^{c}} = (\Id_{\Lie_{1}^{\ash}} \circ \Lambda^{c} \circ \Id_{\qCb^{\ash}}) \cdot (\Delta_{\Lie_{1}^{\ash}} \circ \Delta_{\qCb^{\ash}})$. In order to get the infinitesimal decomposition map, we have to project the result onto $\qCx^{\ash} \circ_{(1)} \qCx^{\ash}$. In the sequel, for three $\Sb$-modules $M$, $N_{1}$ and $N_{2}$, we use the notation $M \circ (N_{1}\, ;\, N_{2})$ for the sub-$\Sb$-module of $M \circ (N_{1} \oplus N_{2})$ linear in $N_{2}$. The infinitesimal decomposition map is given by the composite
\begin{multline*}
\Lie_{1}^{\ash} \circ \qCb^{\ash} \xrightarrow{\Delta_{\Lie_{1}^{\ash}}^{(1)} \circ \Id} (\Lie_{1}^{\ash} \circ_{(1)} \Lie_{1}^{\ash} ) \circ \qCb^{\ash} \cong \Lie_{1}^{\ash} \circ (\qCb^{\ash}\, ;\, \Lie_{1}^{\ash} \circ \qCb^{\ash})\\
\hspace{-3cm} \xrightarrow{\Id \circ (\Id \, ;\, \Id \circ \Delta_{\qCb^{\ash}})} \Lie_{1}^{\ash} \circ (\qCb^{\ash}\, ;\, \Lie_{1}^{\ash} \circ \qCb^{\ash} \circ \qCb^{\ash})\\
\xrightarrow{\Id \circ (\Id \, ;\, \Lambda^{c} \circ \Id)} \Lie_{1}^{\ash} \circ (\qCb^{\ash}\, ;\, \qCb^{\ash} \circ \Lie_{1}^{\ash} \circ \qCb^{\ash}) \cong (\Lie_{1}^{\ash} \circ \qCb^{\ash}) \circ_{(1)} (\Lie_{1}^{\ash} \circ \qCb^{\ash}).
\end{multline*}
A careful computation of this composite gives the infinitesimal decomposition map. Similarly, a careful computation of the composite
$$
\Lie_{1}^{\ash} \circ \qCb^{\ash} \xrightarrow{\Delta_{\Lie_{1}^{\ash}} \circ \Delta_{\qCb^{\ash}}} (\Lie_{1}^{\ash} \circ \Lie_{1}^{\ash}) \circ (\qCb^{\ash} \circ \qCb^{\ash}) \xrightarrow{\Id \circ \Lambda^{c} \circ \Id} (\Lie_{1}^{\ash} \circ \qCb^{\ash}) \circ (\Lie_{1}^{\ash} \circ \qCb^{\ash})
$$
give the full decomposition map.
$\cqfd$
\end{pf}

Finally, we obtain an explicit description of the Kosul dual curved cooperad $\Cx^{\ash}$.

\begin{thm}
The operad $\Cx$ is Koszul and its Koszul dual cooperad \emph{$\Cx^{\ash}$} is isomorphic to the curved cooperad with zero differential \emph{$(\Lie_{1}^{\ash} \circ \qCb^{\ash},\, \Delta_{\Lambda^{c}},\, \theta)$}, where the curvature \emph{$\theta : \Lie_{1}^{\ash} \circ \qCb^{\ash} \rightarrow \Im$} is defined by
$$
\theta \left( (\overline{l}_{n}^{c}\, ;\, \overline{\imath}_{k_{1}}^{c},\, \ldots ,\, \overline{\imath}_{k_{n}}^{c}) \right) = \left\{ \begin{array}{cl}
{\tiny \vcenter{\xymatrix@M=0pt@R=4pt@C=4pt{& \ar@{-}[dddd] &\\ && \\ &&\\ && \\ &&}}} & \textrm{if } (\overline{l}_{n}^{c}\, ;\, \overline{\imath}_{k_{1}}^{c},\, \ldots ,\, \overline{\imath}_{k_{n}}^{c}) = \soqjnrsi ,\\
0 & \textrm{otherwise.}
\end{array} \right.
$$
We get a cofibrant resolution of $\Cx$ \emph{
$$
\Cx_{\infty} := \Omega \Cx^{\ash} \xrightarrow{\sim} \Cx,
$$}
where \emph{$\Omega \Cx^{\ash}$} is the cobar construction on the curved cooperad \emph{$\Cx^{\ash}$}, as defined in \cite{HirshMilles}.
\end{thm}

\begin{rem}
The cobar construction $\Omega \Cx^{\ash}$ is a semi-augmented operad which underlying $\Sb$-module is $\T (s^{-1} \overline{\Cx^{\ash}})$ and which differential is the sum of a term built using the infinitesimal decomposition map of $\Cx^{\ash}$ and a term built using the curvature on $\Cx^{\ash}$.
\end{rem}

\begin{pf}
 The word ``cofibrant'' refers to the model category structure defined by Hinich in \cite{Hinich97, Hinich03}. A second reference is the Appendix A of \cite{MerkulovVallette2}. We refer to Section 4, and more specifically 4.2, of \cite{HirshMilles} for the general theory on curved Koszul duality and for the definition of the Koszul dual curved cooperad. An explicit formula for the decomposition map can be obtained by successive composition of the infinitesimal decomposition map. In Section 4.3 of \cite{HirshMilles}, an operad is said to be Koszul when it has an inhomogeneous presentation satisfying two conditions, called (I) and (II), and when the associated quadratic operad is Koszul. The presentation of $\Cx$ that we give trivially satisfies the two conditions (I) and (II) and the associated quadratic operad $\qCx$ is Koszul by Proposition \ref{Koszul}. It follows that $\Cx$ is Koszul and, from Theorem 4.3.1 in \cite{HirshMilles}, that $\Cx_{\infty}$ is a cofibrant resolution of $\Cx$.
$\cqfd$
\end{pf}

\subsection{Description of the homotopy algebras}

For each application $f : (A,\, d_{A}) \rightarrow (B,\, d_{B})$, we define the differential $\partial (f) := d_{B} \cdot f - (-1)^{|f|} f \cdot d_{A}$. We obtain the following explicit description of homotopy $\Cx$-algebra structure.

\begin{prop}\label{jinfinity}
A $\Cx_{\infty}$-algebra structure on a dg module $(A,\, d_{A})$ is given by a collection of maps, $\{ j_{k_{1},\, \ldots ,\, k_{n}} \}$ for $n \geq 1$ and $k_{1},\, \cdots ,\, k_{n} \geq 0$, with $j_{0} = d_{A}$ and where each $j_{k_{1},\, \ldots ,\, k_{n}}$ is a map $A^{\otimes n} \rightarrow A$ of degree $1-(k_{1}+ \cdots +k_{n})$ such that
\begin{equation}\label{signrule} \tag{i}
j_{k_{1},\, \ldots ,\, k_{n}}(a_{1} \cdots a_{n}) = \varepsilon \times j_{k_{\sigma(1)},\, \ldots ,\, k_{\sigma(n)}}(a_{\sigma(1)} \cdots a_{\sigma(n)}),
\end{equation}
where $\varepsilon$ is given by the Koszul sign rule. (The maps $j_{k_{1},\, \ldots ,\, k_{n}}$ come from the elements $(\overline{l}_{n}^{c}\, ;\, \overline{\imath}_{k_{1}}^{c},\, \ldots ,\, \overline{\imath}_{k_{n}}^{c})$ where the $\overline{\imath}_{k_{i}}^{c}$ is of degree $k_{i}$ ; therefore the sign $\varepsilon$ depends on the $k_{i}$'s and on the $|a_{i}|$'s.) The family of maps $\{ j_{k_{1},\, \ldots ,\, k_{n}} \}$ satisfies the following identities:
\begin{equation} \tag{ii}
\partial (j_{2})(a) = (j_{0} \cdot j_{2} - j_{2} \cdot j_{0})(a) = j_{1}^{2}(a) + a,
\end{equation}
and, when $(n,\, k_{1},\, \ldots ,\, k_{n}) \neq (1,\, 2)$,
\begin{equation} \tag{iii}
\sum_{\substack{p+q=n+1\\ p,\, q \geq 1}} \sum_{\sigma \in Sh_{q,\, p-1}} \sum_{K} \beta^{\sigma}_{k'_{j},\, k''_{j}} \times j_{K_{1}}\left(j_{K_{2}} \left(a_{\sigma (1)} \cdots a_{\sigma (q)} \right) a_{\sigma (q+1)} \cdots a_{\sigma (n)} \right) = 0,
\end{equation}
where $K = \{k'_{\sigma(j)}+k''_{\sigma(j)} = k_{\sigma(j)} \textrm{ ; } k'_{\sigma(j)} = k_{\sigma(j)} \textrm{ for } j > q \textrm{ ; } k' = k'_{\sigma(1)}+ \cdots +k'_{\sigma(q)} \}$, $K_{1} = \{ k',\, k_{\sigma(q+1)},\, \ldots ,\, k_{\sigma(n)} \}$, $K_{2} = \{ k''_{\sigma(1)},\, \ldots ,\, k''_{\sigma(q)} \}$ and $\beta^{\sigma}_{k'_{j},\, k''_{j}} := \alpha^{\sigma}_{k'_{j},\, k''_{j}} \times \varepsilon'$ with $\alpha^{\sigma}_{k'_{j},\, k''_{j}} = \sgn_{k_{1},\, \ldots ,\, k_{n}}\sigma \times \varepsilon_{k'_{j},\, k''_{j}}^{\sigma} \times \alpha_{k'_{\sigma(1)},\, \ldots ,\, k'_{\sigma(q)}}$ is defined in Proposition \ref{Koszul} and $\varepsilon'$ is given by the Koszul rule sign $a_{1} \cdots a_{n} = \varepsilon' a_{\sigma (1)} \cdots a_{\sigma (n)}$.
\end{prop}

\begin{pf}
Since $\Cx_{\infty} = \Omega \Cx^{\ash}$ is a quasi-free operad, a map $j : \Cx_{\infty}(A) \rightarrow A$ of degree 0 is determined by an application $\Cx^{\ash}(A) \rightarrow A$ of degree 1. This is equivalent to a collection of applications $j_{k_{1},\, \ldots ,\, k_{n}} : A^{\otimes n} \rightarrow A$ of degree $1-(k_{1}+ \cdots +k_{n})$, defined by:
$$
j_{k_{1},\, \ldots ,\, k_{n}}(a_{1} \cdots a_{n}) := j ((\overline{l}_{n}^{c}\, ;\, \overline{\imath}_{k_{1}}^{c},\, \ldots ,\, \overline{\imath}_{k_{n}}^{c}) \otimes a_{1} \cdots a_{n}),
$$
and therefore satisfying equations (\ref{signrule}). Moreover, the operad $\Cx_{\infty}$ has a differential $d := d_{0} + d_{2}$, where $d_{0}$ depends on the curvature $\theta$ and $d_{2}$ depends on the infinitesimal decomposition map $\Delta_{\Cx^{\ash}}^{(1)}$. We refer to Section 3.3.5 in \cite{HirshMilles} for more details. Therefore, the fact that $g$ is a dg map gives the additional relations (ii) and (iii) among the applications $j_{k_{1},\, \ldots ,\, k_{n}}$.
$\cqfd$
\end{pf}

\begin{rem}
The maps $j_{k_{1},\, \ldots ,\, k_{n}}$ can equivalently be seen as maps $A[k_{1}] \odot \cdots \odot A[k_{n}] \rightarrow A[1]$.
\end{rem}

\begin{cor}\label{complexstructure}
Let $V$ be a vector space (concentrated in degree 0). A $\Cx_{\infty}$-algebra structure on $V$ is given by a collection of degree 0 maps
$$
\{ \mathring{\jmath}_{n+1} := j_{0,\, \ldots ,\, 0,\, 1} : V^{\odot n} \otimes V \rightarrow V\}_{n\geq 0},
$$
together satisfying the following identities:
\begin{equation} \label{eq1}
\mathring{\jmath}_{1}^{2} = -\Id,
\end{equation}

and for $n>1$,

\begin{equation} \label{eq2}
\sum_{\substack{p+q=n+1\\ p,\, q \geq 1}} \sum_{\substack{\sigma \in Sh_{p-1,\, q}\\ \sigma(n) = n}} \mathring{\jmath}_{p}\left(a_{\sigma(1)} \cdots a_{\sigma (p-1)} \mathring{\jmath}_{q} \left(a_{\sigma (p)} \cdots a_{\sigma (n)} \right) \right) = 0,
\end{equation}

\begin{multline} \label{eq3}
\sum_{\substack{p+q=n+1\\ p \geq 2,\, q \geq 1}} \sum_{\substack{\sigma \in Sh_{p-1,\, q}\\ \sigma(\{p-1,\, n\}) = \{n-1,\, n\}}} \sgn(\sigma_{|\{p-1,\, n\}}) \times \mathring{\jmath}_{p}\left(a_{\sigma(1)} \cdots a_{\sigma (p-2)} \mathring{\jmath}_{q} \left(a_{\sigma(p)} \cdots a_{\sigma (n)} \right) a_{\sigma (p-1)} \right)\\
+ \sum_{\substack{p+q=n+1\\ p \geq 1,\, q \geq 2}} \sum_{\substack{\sigma \in Sh_{p-1,\, q}\\ \sigma(\{n-1,\, n\}) = \{n-1,\, n\}}} \sgn(\sigma_{|\{n-1,\, n\}}) \times \mathring{\jmath}_{p}\left(a_{\sigma(1)} \cdots a_{\sigma (p-1)} \mathring{\jmath}_{q} \left( a_{\sigma (p)} \cdots a_{\sigma (n)} \right) \right) = 0.
%\end{split}
\end{multline}
\end{cor}

\begin{pf}
Be aware of the fact that $V$ is concentrated in degree 0 and that, for instance, $j ((\overline{l}_{3}^{c}\, ;\, \overline{\imath}_{1}^{c},\, \overline{\imath}_{0}^{c} ,\, \overline{\imath}_{0}^{c}) \otimes a_{1} a_{2} a_{3}) = j ((\overline{l}_{3}^{c}\, ;\, \overline{\imath}_{0}^{c},\, \overline{\imath}_{0}^{c} ,\, \overline{\imath}_{1}^{c}) \otimes a_{3} a_{2} a_{1})$. Then, the result is a particular case of Proposition \ref{jinfinity}.
$\cqfd$
\end{pf}

\begin{rems} \leavevmode
\begin{itemize}
\item Pictorially, the map $\mathring{\jmath}_{n}$ corresponds to the $n$-ary tree: ${\tiny \vcenter{\xymatrix@M=0pt@R=4pt@C=5pt{
    1  & \cdots & & n-1 & n\\
    & & & &\\
    & & & & \circ \ar@{-}[u] \\
    \ar@{-}[drr] & \ar@{-}[dr] & \ar@{-}[d] & \ar@{-}[dl] & \ar@{-}[dll] \ar@{-}[u] \\
    & & *{} \ar@{-}[d] & &\\
    & & & &
    }}}$.
\item The maps $\mathring{\jmath}_{n+1}$ can equivalently be seen as degree 0 maps $V^{\odot n} \otimes V[1] \rightarrow V[1]$. The shift of degree on the source space comes from the $\sojnsi$ in the tree presentation of $\mathring{\jmath}_{n}$ and the one on the target space comes from the suspension of the generators $\Cx^{\ash}$ in the cobar construction $\Omega \Cx^{\ash}$.
\end{itemize}
\end{rems}

\subsection{Infinity-morphism of $\Cx_{\infty}$-algebras}

We follow the general theory and we use some notations of Section 6.2 in \cite{HirshMilles}. Let $A$ and $B$ be two $\Cx_{\infty}$-algebras, with structure maps $j^{A}$ and $j^{B}$. For $C = A$ and $B$, we denote by $D_{j^{C}}$ is the curved codifferential on $\Cx^{\ash}(C)$ induced by the differential $d_C$ and by the $\Cx_{\infty}$-algebra structure $j^C$. We remind that a curved codifferential $D_{j^C}$ on $\Cx^{\ash}(C)$ does not necessarily square to 0 and satisfies
$$
D_{j^C}^{2} = \left( \theta \circ \id_{\Cx^{\ash}(C)} \right) \cdot \Delta_{\Cx^{\ash}(C)},
$$
where $\theta$ is the curvature of $\Cx^{\ash}$.
\begin{dei}
An $\infty$-morphism $A \rightsquigarrow B$ of $\Cx_{\infty}$-algebras is a $\Cx^{\ash}$-coalgebra map
$$
f : \left( \Cx^{\ash}(A),\, D_{j^{A}} \right) \rightarrow \left( \Cx^{\ash}(B),\, D_{j^{B}} \right),
$$
commuting with the curved codifferentials.
\end{dei}

\begin{prop}\label{infinitymorphism}
Let $(A,\, j^{A})$, $(B,\, j^{B})$ be two $\Cx_{\infty}$-algebras. An $\infty$-morphism between $A$ and $B$ is a collection of maps,
$$
\left\{ f_{k_{1},\, \ldots ,k_{n}} : A^{\otimes n} \rightarrow B \right\}_{k_{1},\, \ldots ,\, k_{n} \geq 0,\, n \geq 1} \textrm{ of degree } - (k_{1}+ \cdots + k_{n}),
$$
satisfying:
\begin{multline*}
\partial (f_{k_{1},\, \ldots ,k_{n}}) =\\
\sum_{\substack{p+q=n+1\\ p,\, q \geq 1}} \sum_{\sigma \in Sh_{q,\, p-1}} \sum_{\substack{k'_{\sigma(j)}+k''_{\sigma(j)} = k_{\sigma(j)}\\ k'_{\sigma(j)} = k_{\sigma(j)} \textrm{ for } j > q\\ k' = k'_{\sigma(1)}+ \cdots +k'_{\sigma(q)}}} \alpha^{\sigma,\, \kappa}_{k'_{j},\, k''_{j}} \times f_{k',\, k_{\sigma(q+1)},\, \ldots ,\, k_{\sigma(n)}} \cdot \left(j^A_{k''_{\sigma(1)},\, \ldots ,\, k''_{\sigma(q)}},\, \Id_A ,\, \ldots \right)^{\sigma^{-1}}\\
- \sum_{\substack{\{q_{1},\, \ldots ,\, q_{p} \}\\ q_{1}+\cdots +q_{p}= n}} \sum_{\sigma \in Sh_{q_{1},\, \ldots ,\, q_{p}}} \frac{1}{N_{\underline{q}}} \sum_{k'_{j}+k''_{j} = k_{j}} \beta^{\sigma,\, \gamma}_{k'_{j},\, k''_{j}} \times j^B_{l'_1,\, \ldots ,\, l'_p} \cdot \left( f_{k''_{\sigma(1)},\, \ldots ,\, k''_{\sigma(q_1)}},\, f_{k''_{\sigma(q_1+1)},\, \ldots ,\, k''_{\sigma(q_1+q_2)}} ,\, \ldots \right)^{\sigma^{-1}}
\end{multline*}
where $\alpha^{\sigma,\, \kappa}_{k'_{j},\, k''_{j}} := (-1)^{k'+k_{\sigma(q+1)}+\cdots +k_{\sigma(n)}} \times \alpha^{\sigma}_{k'_{j},\, k''_{j}}$ and $\beta^{\sigma,\, \gamma}_{k'_{j},\, k''_{j}} := (-1)^{\left(-\sum k''_j \right) \left(\sum k'_j \right)} \times \beta^{\sigma}_{k'_{j},\, k''_{j}}$ and where $N_{\underline{q}} := \prod_{i=1}^{m} \max(n_{i},\, 1)$ for $\underline{q} = \{ q_{1},\, \ldots ,\, q_{n}\} = \{ \underbrace{1,\, \ldots ,\, 1}_{n_{1} \textrm{ times}},\, 2,\, \ldots ,\, \underbrace{m,\, \ldots ,\, m}_{n_{m} \textrm{ times}} \}$.
\end{prop}

\begin{pf}
A $\Cx^{\ash}$-coalgebra map $f : \Cx^{\ash}(A) \rightarrow \Cx^{\ash}(B)$ (of degree 0) is characterized by its corestriction to $B$, that is $f$ is determined by a collection of maps $f_{k_{1},\, \ldots ,k_{n}} : A^{\otimes n} \rightarrow B$ of degree $-(k_1+ \cdots + k_n)$.The fact that $f$ commutes with the curved codifferentials is equivalent to the following commutative diagram
$$
\xymatrix@C=50pt{\Cx^{\ash} (A) \ar[r]^{\Delta \circ \Id_A} \ar[d]_{d_1+d_2} & \Cx^{\ash} \circ \Cx^{\ash} (A) \ar[r]^{\Id_{\Cx^{\ash}} \circ f} & \Cx^{\ash}(B) \ar[d]^{d_B+\bar{d}_2}\\
\Cx^{\ash} (A) \ar[rr]_f && B,}
$$
where $d_1$, induced by the differential $d_A$ or $d_B$, and $d_2$, induced by the algebra structure $j^A$ or $j^B$, are defined in Section 5.2.3 of \cite{HirshMilles} for example, and $\bar{d}_2$ is the projection of $d_2$ on $B$. Making this diagram explicit gives exactly the conditions concerning the maps $f_{k_1,\, \ldots ,\, k_n}$.
$\cqfd$
\end{pf}

\begin{cor}\label{corinfinitymorph}
Let $V$ and $W$ be two vector spaces (concentrated in degree 0) endowed with $\Cx_{\infty}$-algebra structures $\mathring{\jmath}^V$ and $\mathring{\jmath}^W$. An $\infty$-morphism between $V$ and $W$ is given by a collection of degree 0 maps
$$
\{ \mathring{f}_n := f_{0,\, \ldots ,\, 0} : V^{\otimes n} \rightarrow W\}_{n\geq 1},
$$
satisfying, on elements $v_1$, $\ldots$, $v_n$ in $V$,
\begin{multline}\label{infinitymorph}
\sum_{\substack{p+q=n+1\\ p,\, q \geq 1}} \sum_{\substack{\sigma \in Sh_{q,\, p-1}\\ \sigma(q)=n}} \mathring{f}_p \cdot \left(\mathring{\jmath}^V_q (v_{\sigma(1)} \cdots v_{\sigma(q)}) v_{\sigma(q+1)} \cdots v_{\sigma(n)} \right) =\\
\sum_{\substack{\{q_{1},\, \ldots ,\, q_{p-1} \},\, q_p \\ q_{1}+\cdots +q_{p}= n}} \frac{N^{q_p}_{\underline{q}}}{N_{\underline{q}}} \sum_{\substack{\sigma \in Sh_{q_{1},\, \ldots ,\, q_{p}}\\ \sigma(n)=n}} \mathring{\jmath}^W_p \cdot \left( \mathring{f}_{q_1}(v_{\sigma(1)} \cdots v_{\sigma(q_1)}) \mathring{f}_{q_2}(v_{\sigma(q_1+1)} \cdots v_{\sigma(q_1+q_2)}) \ldots \right),
\end{multline}
where $N^{q_p}_{\underline{q}} := \mathrm{Card} \left\{ q \in \underline{q} = \{q_1,\, \ldots ,\, q_p\} \textrm{ s.t. } q= q_p \right\}$ and where $N_{\underline{q}}$ is defined in Proposition \ref{infinitymorphism}.
\end{cor}

\begin{pf}
The result is a particular case of Proposition \ref{infinitymorphism}. The factors $N_{\underline{q}}^{q_{p}}$ comes from the fact that we fix $\sigma(n) =n$.
$\cqfd$
\end{pf}

\begin{rem}
Pictorially, we get
$$
\sum {\tiny \vcenter{\xymatrix@M=2pt@R=4pt@C=5pt{
    & & & & \\
    & & *{\circ} \ar@{-}[u] & & & \\
    & & *{} \ar@{-}[u] & & & \\
    & \mathring{\jmath}_2^V \ar@{-}[d] \ar@{-}[ul] \ar@{-}[ur] & & & & \\
    & *{} \ar@{-}[dr] & \ar@{-}[d] & \ar@{-}[dl] & \ar@{-}[dll] & \ar@{-}[dlll] \\
    & & \mathring{f}_5 \ar@{-}[d] & & &\\
    & & & & &
    }}} = \sum {\tiny \vcenter{\xymatrix@M=2pt@R=4pt@C=5pt{
    & & & & &\\
    \ar@{-}[dr] & \ar@{-}[d] & \ar@{-}[dl] & \ar@{-}[dr] & & \ar@{-}[dl] \\
    & \mathring{f}_3 & & & \mathring{f}_2 & \\
    & *{} \ar@{-}[d] \ar@{-}[u] & & \ar@{-}[dd] & *{\circ} \ar@{-}[u] & \\
    & *{} \ar@{-}[drr] & & & *{} \ar@{-}[dl] \ar@{-}[u] & \\
    & & & \mathring{\jmath}_3^W \ar@{-}[d] & &\\
    & & & & &
    }}}.
$$
\end{rem}

\section{Formal complex structures are homotopy algebras}\label{section3}

In this section, we give a geometric interpretation of the $\Cx_{\infty}$-algebras. We prove that a $\Cx_{\infty}$-algebra structure on a vector space $V \cong \Rb^m$ is precisely a complex structure on the formal manifold $\V_{\mathrm{for}} \cong \left(\{ 0\},\, S(\Rb^m)^*\right)$. We also show that $\infty$-morphisms of $\Cx_{\infty}$-algebras correspond to holomorphic maps of formal complex structures.\\

In order to prove this result, we fix the following notations. Let $\left\{ e_{a}^{V} \right\}_{a}$, resp. $\left\{ e_{a}^{W} \right\}_{a}$, be a basis of $V$, resp. $W$. We denote by $j^{b}_{\alpha_{1},\, \ldots ,\, \alpha_{n},\, a}$ the coefficients of $\jo_{n+1} : V^{\odot n} \otimes V \rightarrow V$ on the basis
$$
\jo_{n+1}(e_{\alpha_{1}}^{V},\, \ldots ,\, e_{\alpha_{n}}^{V}\, ;\, e_{a}^{V}) = j^{b}_{\alpha_{1},\, \ldots ,\, \alpha_{n},\, a} e_{b}^{V},
$$
and by $f^{b}_{\alpha_{1},\, \ldots ,\, \alpha_{n}}$ the coefficients of $\fo_{n} : V^{\odot n} \rightarrow W$ on the bases
$$
\fo_{n}(e_{\alpha_{1}}^{V},\, \ldots ,\, e_{\alpha_{n}}^{V}) = f^{b}_{\alpha_{1},\, \ldots ,\, \alpha_{n}} e_{b}^{W}.
$$
Then, to any family of maps $\{ \jo_{n+1} : V^{\odot n} \otimes V \rightarrow V \}_{n \geq 0}$, we associate the endomorphism $J$ of the tangent sheaf $\T_{\V}[1]$ defined by:
$$
J := \sum_{n \geq 0} j_{\alpha_{1},\, \ldots ,\, \alpha_{n},\, a}^{b} \frac{1}{\underline{\alpha}!} t^{\alpha_{1}} \cdots t^{\alpha_{n}} \gamma^{a} \partial_{b} = J_{a}^{b}(t) \gamma^{a} \partial_{b},
$$
where $\underline{\alpha} = \{ \alpha_{1},\, \ldots ,\, \alpha_{n}\} = \{ \underbrace{1,\, \ldots ,\, 1}_{n_{1} \textrm{ times}},\, 2,\, \ldots ,\, \underbrace{p,\, \ldots ,\, p}_{n_{p} \textrm{ times}} \}$ and $\underline{\alpha}! := \prod_{i = 1}^{p} n_{i} !$. We denote moreover $N_{\underline{\alpha}} := \prod_{i=1}^{p} n_{i}$. Let fix a subscript $c$ in $\underline{\alpha}$. Since $J_{a}^{b}$ is symmetric in the $t^{\alpha_{i}}$'s, we can assume that the product $t^{\alpha_{1}} \cdots t^{\alpha_{n}}$ satisfies $\alpha_{1} = c$. We have the following computation:
$$
\partial_{c} J_{a}^{b} (t) = \sum_{\underline{\alpha} \textrm{ s.t. } c \in \underline{\alpha}} j_{c,\, \alpha_{2},\, \ldots ,\, \alpha_{n},\, a}^{b} \frac{1}{\{ \alpha_{2},\, \ldots ,\, \alpha_{n}\}!} t^{\alpha_{2}} \cdots t^{\alpha_{n}}.
$$
Similarly, to any family of maps $\{ \fo_{n} : V^{\odot n} \rightarrow W \}_{n\geq 0}$, we associate the map of formal manifolds $F : \V_{\mathrm{for}} = \left(\{\textrm{point} \},\, \O_{\V_{\mathrm{for}}}\right) \rightarrow \W_{\mathrm{for}} =  \left(\{\textrm{point} \},\, \O_{\W_{\mathrm{for}}}\right)$ defined by:
$$
F := \sum_{n \geq 1} f_{\alpha_{1},\, \ldots ,\, \alpha_{n}}^{b} \frac{1}{\underline{\alpha}!} t^{\alpha_{1}} \cdots t^{\alpha_{n}} e_{b}^{W} = F_{a}^{b}(t) e_{b}^{W}.
$$

We denote by $\Vect$ the category of finite-dimensional vector spaces (concentrated in degree 0).

\begin{thm}\label{equiv1}
There is an equivalence of categories
$$
\begin{array}{ccc}
\left\{ \begin{gathered} \Cx_{\infty}\textrm{-algebra structures on } \Vect \\ \textrm{with } \infty \textrm{-morphisms} \end{gathered} \right\} & \xrightarrow{\cong} & \left\{ \begin{gathered} \textrm{Complex structures on formal manifolds}\\ \V_{\mathrm{for}} = \left(\{\textrm{point} \},\, S(V)^*\right) \textrm{ with } V \in \Vect \\ \textrm{with holomorphic maps} \end{gathered} \right\}\\
&&\\
{\displaystyle \left\{ \jo_{n+1} : V^{\odot n} \otimes V \rightarrow V \right\}_{n\geq 0}} & \mapsto & {\displaystyle J = \sum_{n \geq 0} j_{\alpha_{1},\, \ldots ,\, \alpha_{n},\, a}^{b} \frac{1}{\underline{\alpha}!} t^{\alpha_{1}} \cdots t^{\alpha_{n}} \gamma^{a} \partial_{b}}\\
&&\\
{\displaystyle \left\{ \fo_n : V^{\otimes n} \rightarrow W \right\}_{n\geq 0}} & \mapsto & {\displaystyle F := \sum_{n \geq 1} f_{\alpha_{1},\, \ldots ,\, \alpha_{n}}^{b} \frac{1}{\underline{\alpha}!} t^{\alpha_{1}} \cdots t^{\alpha_{n}} e_{b}^{W}.}
\end{array}
$$
\end{thm}

\begin{pf}
The composition of $J$ with itself gives the following computation
\begin{equation*}
\begin{split}
J^{2} & = (J_{a}^{b}(t) \gamma^{a} \partial_{b}) \cdot (J_{a'}^{b'}(t) \gamma^{a'} \partial_{b'}) = J_{a}^{b}(t) J_{a'}^{a}(t) \gamma^{a'} \partial_{b}\\
& = \left( \sum_{K \sqcup L = \{ \alpha_{1},\, \ldots ,\, \alpha_{n} \}} j_{K,\, a'}^{a} \frac{1}{K!} j_{L,\, a}^{b} \frac{1}{L!} t^{\alpha_{1}} \cdots t^{\alpha_{n}} \right) \gamma^{a'} \partial_{b}\\
& = \sum_{K \sqcup L = \{ \alpha_{1},\, \ldots ,\, \alpha_{n} \}} \left( j_{K,\, a'}^{a}\ j_{L,\, a}^{b} \frac{K \sqcup L\ !}{K! \cdot L!} \right) \cdot \frac{1}{K \sqcup L\ !} t^{\alpha_{1}} \cdots t^{\alpha_{n}} \gamma^{a'} \partial_{b}.
\end{split}
\end{equation*}
The coefficient in the parentheses is the coefficient of $e_{b}$ in the left-hand side of Equations (\ref{eq1}) and (\ref{eq2}) in Corollary \ref{complexstructure}, with $n+1$ instead of $n$, $a_{n+1} := e_{a'}$ and $a_{i} := e_{\alpha_{i}}$ for $i \leq n$, since the number of shuffles is given by the binomial coefficient. The equation $J^{2} = -\Id$ is therefore equivalent to Equations (\ref{eq1}) and (\ref{eq2}).\\

Similarly, we have
\begin{multline*}
\left( J_{a}^{b}(t) \partial_{b} J_{a'}^{d}(t) - J_{a'}^{b}(t) \partial_{b} J_{a}^{d}(t) \right) \gamma^{a} \gamma^{a'} \partial_{d} =\\
= \left( \sum_{K \sqcup L = \{ \alpha_{1},\, \ldots ,\, \alpha_{n} \}} \left( j_{b,\, K,\, a'}^{d} \frac{1}{K!} j_{L,\, a}^{b} \frac{1}{L!} - j_{b,\, K,\, a}^{d} \frac{1}{K!} j_{L,\, a'}^{b} \frac{1}{L!} \right) t^{\alpha_{1}} \cdots t^{\alpha_{n}}\right) \gamma^{a} \gamma^{a'} \partial_{d}\\
= \sum_{K \sqcup L = \{ \alpha_{1},\, \ldots ,\, \alpha_{n} \}} \left( \left( j_{b,\, K,\, a'}^{d}\ j_{L,\, a}^{b} - j_{b,\, K,\, a}^{d}\ j_{L,\, a'}^{b} \right) \frac{K \sqcup L\ !}{K! \cdot L!} \right) \cdot \frac{1}{K \sqcup L\ !} t^{\alpha_{1}} \cdots t^{\alpha_{n}} \gamma^{a} \gamma^{a'} \partial_{d}.
\end{multline*}
The coefficient in the parentheses is the opposite of the coefficient of $e_{b}$ in the first term of the left-hand side of Equation (\ref{eq3}) in Corollary \ref{complexstructure}, with $n+2$ instead of $n$, $,a_{n+1} := e_{a'}$, $a_{n+2} := e_{a}$ and $a_{i} := e_{\alpha_{i}}$ for $i \leq n$. Moreover, we have
\begin{multline*}
\left( J_{b}^{d}(t) \partial_{a} J_{a'}^{b}(t) - J_{b}^{d}(t) \partial_{a'} J_{a}^{b}(t) \right) \gamma^{a} \gamma^{a'} \partial_{d} =\\
= \left( \sum_{K \sqcup L = \{ \alpha_{1},\, \ldots ,\, \alpha_{n} \}} \left( j_{a,\, K,\, a'}^{b} j_{L,\, b}^{d} - j_{a',\, K,\, a}^{b} j_{L,\, b}^{d} \right) \frac{1}{K!} \frac{1}{L!} t^{\alpha_{1}} \cdots t^{\alpha_{n}}\right) \gamma^{a} \gamma^{a'} \partial_{d}\\
= \sum_{K \sqcup L = \{ \alpha_{1},\, \ldots ,\, \alpha_{n} \}} \left( \left( j_{a,\, K,\, a'}^{b}\ j_{L,\, b}^{d} - j_{a',\, K,\, a}^{b}\ j_{L,\, b}^{d} \right) \frac{K \sqcup L\ !}{K! \cdot L!} \right) \cdot \frac{1}{K \sqcup L\ !} t^{\alpha_{1}} \cdots t^{\alpha_{n}} \gamma^{a} \gamma^{a'} \partial_{d}.
\end{multline*}
The coefficient in the parentheses is the opposite of the coefficient of $e_{b}$ in the second term of the left-hand side of Equation (\ref{eq3}) in Corollary \ref{complexstructure}, with $n+2$ instead of $n$, $a_{n+1} := e_{a}$, $a_{n+2} := e_{a'}$ and $a_{i} := e_{\alpha_{i}}$ for $i \leq n$. By means of equality (\ref{bracket}), we obtain that Equation (\ref{eq3}) is equivalent to the integrability condition and therefore, Equations (\ref{rel}) are equivalent to Equations (\ref{eq1}), (\ref{eq2}) and (\ref{eq3}).\\

It remains to check the equivalence on morphisms. We have
\begin{equation*}
\begin{split}
\partial_{b} F^{d}(t) J_{a}^{b}(t) \gamma^{a} \partial_{d} & = \left( \sum_{K \sqcup L = \{ \alpha_{1},\, \ldots ,\, \alpha_{n} \}} f_{b,\, K}^{d} \frac{1}{K!} j_{L,\, a}^{b} \frac{1}{L!} t^{\alpha_{1}} \cdots t^{\alpha_{n}}\right) \gamma^{a} \partial_{d}\\
& = \sum_{K \sqcup L = \{ \alpha_{1},\, \ldots ,\, \alpha_{n} \}} \left( f_{b,\, K}^{d}\ j_{L,\, a}^{b} \frac{K \sqcup L\ !}{K! \cdot L!} \right) \cdot \frac{1}{K \sqcup L\ !} t^{\alpha_{1}} \cdots t^{\alpha_{n}} \gamma^{a} \partial_{d}. 
\end{split}
\end{equation*}
The coefficient in the parentheses is equal to the coefficient of $e_{d}^{W}$ in the left-hand side of Equation (\ref{infinitymorph}) in Corollary \ref{corinfinitymorph} with $n+1$ instead of $n$, $v_{i} := e_{\alpha_{i}}^{V}$ for $i<q$, $v_{q} := e_{a}^{V}$ and $v_{i} := e_{\alpha_{i-1}}^{V}$ for $i>q$. And finally, we compute
\begin{multline*}
{J'}_{b}^{d}(F(t))  \partial_{a} F^{b}(t) \gamma^{a} \partial_{d}\\
\shoveleft{= \left( \sum_{\substack{K \sqcup L = \{ \alpha_{1},\, \ldots ,\, \alpha_{n} \}\\ K = K_{1} \sqcup \cdots \sqcup K_{p}\\ B = \{b_{1},\, \ldots ,\, b_{p} \}}} {j'}_{b_{1},\, \ldots ,\, b_{p},\, b}^{d} \frac{1}{B !} f_{K_{1}}^{b_{1}} \frac{1}{K_{1}!} \cdots f_{K_{p}}^{b_{p}} \frac{1}{K_{p}!} f_{a,\, L}^{b} \frac{1}{L !} t^{\alpha_{1}} \cdots t^{\alpha_{n}}\right) \gamma^{a} \partial_{d}}\\
= \sum_{\substack{K \sqcup L = \{ \alpha_{1},\, \ldots ,\, \alpha_{n} \}\\ K = K_{1} \sqcup \cdots \sqcup K_{p}\\ B = \{b_{1},\, \ldots ,\, b_{p} \}}} \left( {j'}_{b_{1},\, \ldots ,\, b_{p},\, b}^{d} f_{K_{1}}^{b_{1}} \cdots f_{K_{p}}^{b_{p}} \frac{1}{B !} \cdot f_{a,\, L}^{b} \frac{K \sqcup L\ !}{K_{1}! \cdots K_{p}! \cdot L!} \right) \cdot \frac{1}{K \sqcup L\ !} t^{\alpha_{1}} \cdots t^{\alpha_{n}} \gamma^{a} \partial_{d}.
\end{multline*}
The coefficient in the parentheses is equal to the coefficient of $e_{d}^{W}$ in the right-hand side of Equation (\ref{infinitymorph}) in Corollary \ref{corinfinitymorph} with $n+1$ instead of $n$, $v_{i} := e_{\alpha_{i}}^{V}$ for $i \leq n$ and $v_{n+1} := e_{a}^{V}$. The coefficient $\frac{N^{q_p}_{\underline{q}}}{N_{\underline{q}}}$ does not appear here since each term appears already only once, in comparison to the formula for the full decomposition map of $\Lie_{1}^{\ash}$ for instance. Therefore, Equation (\ref{holom}) is equivalent to Equation (\ref{infinitymorph}) in Corollary \ref{corinfinitymorph} and this concludes the proof.
$\cqfd$
\end{pf}

\section{Globalisation}\label{section4}

We devote this last section to the globalisation of the results of the previous section. A complex manifold can be described locally in terms of coordinates and this is the reason why we introduce the space of all coordinates systems on a manifold. We build a fiber bundle endowed with a connection in order to characterize the smooth complex structures on a manifold as certain families of $\Cx_{\infty}$-algebras. We propose a similar result for morphisms.

\subsection{Coordinate space and connexion}\label{coorandconn}

We use here the work of \cite{BernshteinRozenfeld} and follow the ideas of \cite{CattaneoFelderTomassini}. Let $M$ be a paracompact smooth $m$-manifold. We denote by $\mathcal{M}$ the locally ringed space $(M,\, \mathcal{C}^{\infty}_{M})$, for example $\R^m = (\Rb^m,\, \mathcal{C}^{\infty}_{\Rb^m})$, and by $\mathcal{R}^{m}_{\mathrm{for}}$ the formal pointed manifold associated to the vector space $\Rb^{m}$. We consider the manifold $M^{\textrm{coor}}$ given by
$$
M^{\textrm{coor}} := \left\{ \begin{array}{l}
(x,\, \varphi) \textrm{ such that } x\in M \textrm{ and } \varphi : \mathcal{R}^{m}_{\mathrm{for}} \rightarrow (\mathcal{M},\, x)\\
\textrm{ is a pointed immersion of locally ringed spaces} \end{array} \right\}.
$$
It can be thought as the space of all local coordinate system on $M$. In this section, we construct a connexion on the trivial vector bundle $\tE$ over $M^{\textrm{coor}}$,
$$
\tE := M^{\textrm{coor}} \times \mathrm{M}_{m}\left(\Rb \llbracket t^1,\, \ldots ,\, t^m \rrbracket \right) \xrightarrow{pr_{1}} M^{\textrm{coor}},
$$
whose fibers over each point are $m \times m$-matrices with coefficients in the formal power series ring $\Rb \llbracket t^1,\, \ldots ,\, t^m \rrbracket$.\\

 Let $G^{0}_m := \mathrm{Aut}\left(\mathcal{R}^{m}_{\mathrm{for}} \right)$ be the (pro-Lie) group of pointed formal diffeomorphism of $\Rb^{m}$. The fiber bundle $\Mco$ over $M$ is a principal $G^{0}_m$-bundle, whose left-action is given by
$$
\rho : G^{0}_m \times \Mco \rightarrow \Mco,\ (g,\, (x,\, \varphi)) \mapsto (x,\, \varphi \cdot g^{-1}).
$$
The derivative of $\rho_{(x,\, \varphi)} : G^{0}_m \rightarrow \Mco$ provides a Lie algebra morphism $T_{\Id}G^{0}_m \rightarrow \chi(\Mco)$, where $\chi(\Mco)$ is the space of vector fields on $\Mco$. We denote by $J^{\infty}(M,\, x)$ the space of infinite jet of vector fields on $M$ at $x$. For example, we have
$$
W_m := J^{\infty}(\Rb^{m},\, 0) = \left\{ v^{b}(t) \partial_{b}\ |\ v^{b}(t) \in \Rb \llbracket t^{1},\, \ldots ,\, t^{m} \rrbracket \right\}.
$$
It is the space of vector fields $\chi \left(\mathcal{R}^{m}_{\mathrm{for}} \right)$ on $\R^{m}_{\mathrm{for}}$. We call them \emph{formal vector fields} (on $\Rb^m$ at $0$). We denote by $W^{0}_m$ the subspace of formal vector fields vanishing at the point $0$. It corresponds via the following isomorphism $s_{(0,\, \Id)}$ to the tangent space of $G^{0}_m$ at $\mathrm{Id}$.

\begin{thm}[Theorem 4.1 of \cite{BernshteinRozenfeld}]
We have the sequence of linear isomorphisms
$$
T_{(0,\, \Id)}(\Rb^{m})^{\mathrm{coor}} \xleftarrow{s_{(0,\, \Id)}} W_m \xleftarrow{\alpha_{(x,\, \varphi)}} J^{\infty}(M,\, x) \xrightarrow{s_{(x,\, \varphi)}} T_{(x,\, \varphi)} M^{\mathrm{coor}},
$$
where the maps $s_{(x,\, \varphi)}$ (denoted $\sigma_x$ in \cite{BernshteinRozenfeld}) can be seen as lifting homomorphisms and the map $\alpha_{(x,\, \varphi)}$ is characterized by the map $\varphi$. (It sends the infinite jet of a vector field to its pushforward by $\varphi^{-1}$.)

The composition $\beta_{(x,\, \varphi)} := - s_{(x,\, \varphi)} \cdot \alpha_{(x,\, \varphi)}^{-1}$ defines a morphism of Lie algebras $W_m \rightarrow \chi(\Mco)$ which extends the previous morphism $W^{0}_m \rightarrow \chi(\Mco)$. Taking its inverse at each point, we get a $W_m$-valued differential form $\omega \in \Omega^1(\Mco,\, W_m)$ which is invariant under the action of diffeomorphisms of $\Mco$ induced by diffeomorphisms of $M$.
\end{thm}

Let $(x,\, \varphi) \in \Mco$ and $\tilde{\varphi} : \R^m \rightarrow \M$ be a local diffeomorphism which gives $\varphi$ when restricted to $\R^m_{\mathrm{for}}$. We define the map
$$
\tilde{\varphi}^{\mathrm{coor}} : (\Rb^m)^{\mathrm{coor}} \rightarrow \Mco,\ (t,\, \tau_t) \mapsto (\tilde{\varphi}(t),\, \tilde{\varphi} \cdot \tau_t).
$$
The derivative of this map at the point $(0,\, \Id)$ does not depend on the chosen $\tilde{\varphi}$. We therefore denote it by $d \varphi^{\coor}(0,\, \Id) : T_{(0,\, \Id)} (\Rb^m)^{\coor} \rightarrow T_{(x,\, \varphi)} \Mco$. By construction, we have the commutative diagram
\begin{equation}\label{diagramcomm}
\begin{gathered}
\xymatrix@C=60pt@R=16pt{W_m \ar[d]_{s_{(0,\, \Id)}} & J^{\infty}(M,\, x) \ar[l]_{\alpha_{(x,\, \varphi)}} \ar[d]^{s_{(x,\, \varphi)}}\\
T_{(0,\, \Id)}(\Rb^{m})^{\mathrm{coor}} \ar[r]_{d \varphi^{\coor}(0,\, \Id)} & T_{(x,\, \varphi)} M^{\mathrm{coor}}.}
\end{gathered}
\end{equation}

The differential $d_t g$ of a diffeomorphism $g$ in $G^{0}_m$ is an $\O_{\R^{m}_{\mathrm{for}}}$-linear map $T_{\R^{m}_{\mathrm{for}}} \rightarrow T_{\R^{m}_{\mathrm{for}}}$. A matrix $A(t) \in \Mrm_{m}\left(\Rb \llbracket t^1,\, \ldots ,\, t^m \rrbracket \right)$ can also be seen as an  $\O_{\R^{m}_{\mathrm{for}}}$-linear map $T_{\R^{m}_{\mathrm{for}}} \rightarrow T_{\R^{m}_{\mathrm{for}}}$.

\begin{dei}\label{rem1}
For $g \in G^{0}_m$ and $A(t) \in \Mrm_{m}\left(\Rb \llbracket t^1,\, \ldots ,\, t^m \rrbracket \right)$, we define the left-action of $G^{0}_m$ on $\Mrm_{m}\left(\Rb \llbracket t^1,\, \ldots ,\, t^m \rrbracket \right)$ by
$$
g \ast A(t) := d_t g (t \cdot g^{-1}) \cdot A(t \cdot g^{-1}) \cdot (d_t g (t \cdot g^{-1}))^{-1},
$$
where $t \cdot g^{-1} = (t^1,\, \ldots ,\, t^m) \cdot g^{-1} := (t^1 \cdot g^{-1},\, \ldots ,\, t^m \cdot g^{-1})$ is the composite of the $t^k$'s seen as an application $\R^m_{\mathrm{for}} \rightarrow \R_{\mathrm{for}}$ and $g^{-1} \in \mathrm{Aut}(\R^m_{\mathrm{for}})$. The derivative of the Lie action provides a linear map $\ast'$ and therefore, we get a linear map $\tilde{\ast}'$
$$
W^0_m \xrightarrow{s_{(0,\, \Id)}} T_{\Id} G^{0}_m \xrightarrow{\ast'} \chi \left(\Mrm_{m}\left(\Rb \llbracket t^1,\, \ldots ,\, t^m \rrbracket \right)\right).
$$
\end{dei}

\begin{rem}
Explicitly, for any element $A(t)$ in $\Mrm_{m}\left(\Rb \llbracket t^1,\, \ldots ,\, t^m \rrbracket \right)$, the derivative of the map $g \mapsto A(t \cdot g^{-1})$ associates to any $\xi \in W^0_m$ the matrix $(d_{t} A)(t \cdot g^{-1}) \cdot d_g (g^{-1}) (\xi)$ and the derivative of the map $g \mapsto d_t g \cdot A \cdot (d_t g)^{-1}$ associates the matrix $d_g (d_t g)(\xi) \cdot A \cdot (d_t g)^{-1} + d_t g \cdot A \cdot d_g ((d_t g)^{-1}) (\xi)$.
\end{rem}

A matrix in $\Mrm_{m}\left(\Rb \llbracket t^1,\, \ldots ,\, t^m \rrbracket \right)$ is an endomorphism $\R^m_{\mathrm{for}} \rightarrow $  A constant formal vector field $v^{b}(0) \partial_{b}$ acts naturally on $\Mrm_{m}\left(\Rb \llbracket t^1,\, \ldots ,\, t^m \rrbracket \right)$ by differentiation of all the coefficients by $v^{b}(0) \partial_{b}$. We can therefore extend the map $\tilde{\ast}'$ to
$$
W_m \xrightarrow{\tilde{\ast}'} \chi \left(\Mrm_{m}\left(\Rb \llbracket t^1,\, \ldots ,\, t^m \rrbracket \right)\right).
$$

Let $U$ be a contractible open subset of $M$. There exist sections $U \rightarrow U^{\coor}$ and we fix $\psi$ such a section. To a tangent vector $\xi_u \in T_u U$, we associate the vector field $\hat{\xi}_u(\psi) := (\psi^* \omega)(\xi_u) \in W_m$, where $\psi^* \omega : TU \rightarrow W_m $ is the pullback of the form $\omega$ by $\psi$.

\begin{prop}\label{conn}
We define a connection on the trivial bundle
$$
E_U := U \times \mathrm{M}_{m}\left(\Rb \llbracket t^1,\, \ldots ,\, t^m \rrbracket \right) \xrightarrow{pr_{1}} U
$$
as follows:
$$
\begin{array}{lccl}
\nabla_{U}(\psi) : & \Gamma(E_U) & \rightarrow & \Omega_{U}^{1} \otimes_{\C^{\infty}(U)} \Gamma(E_U)\\
& \sigma & \mapsto & d_{dR} \sigma + \omega_{U}(\sigma),
\end{array}
$$
where $d_{dR}$ is the de Rham differential and where the application $\omega_U$ is defined, for any vector field $\xi \in \chi_U$, by $\omega_{U}(\sigma)(\xi)(u) := \tilde{\ast}'(\hat{\xi}_{u}(\psi))(\sigma(u))$.
\end{prop}

\begin{pf}
The map $\omega_U (\sigma)(-) : \chi_U \rightarrow \Gamma (E_U)$ is defined point-wise by linear maps, hence it is $\C^{\infty}(U^{\coor})$-linear and $\nabla_{U}(\psi)$ is well-defined. The map $\omega_{U}(-)(\xi)$ is $\C^{\infty}(U^{\coor})$-linear since the left-action $\ast$ is linear in the second variable. It follows that $\nabla_{U}(\psi)$ is a connection.
$\cqfd$
\end{pf}

\begin{lem}\label{comp}
Let $\psi,\, \psi' : U \rightarrow U^{\coor}$ be two sections of $U^{\coor}$ such that $\psi'(u)(t) = \psi(u) (t \cdot g(u)^{-1})$ for some smooth map $g : U \rightarrow G^0_m$, and let $\nabla_U = \nabla_U(\psi)$ and $\nabla'_U = \nabla_U(\psi')$ the corresponding connections. Then, for any section $\sigma \in \Gamma(E_U)$, we have
$$
\nabla'_U (g \ast \sigma) = g \ast \nabla_U \sigma,
$$
where $G^0_m$ acts on $\Omega_{U}^{1} \otimes_{\C^{\infty}(U)} \Gamma(E_U)$ by means of its action on the second factor.
\end{lem}

\begin{pf}
Let $\xi \in \chi_U$ be a vector field. First, we have\\
$d_{dR}(g \ast \sigma)(\xi) = d_{dR}\left(d_t g (t \cdot g^{-1}) \cdot \sigma(t \cdot g^{-1}) \cdot (d_t g (t \cdot g^{-1}))^{-1}\right)(\xi)$
$$
\begin{array}{clc}
= & \left((d_x d_t g) (t \cdot g^{-1})(\xi) + (d_t d_t g)(t \cdot g^{-1}) \cdot d_x (g^{-1})(\xi)\right) \cdot \sigma(t \cdot g^{-1}) \cdot (d_t g (t \cdot g^{-1}))^{-1} & +\\
& d_t g (t \cdot g^{-1}) \cdot \left( d_x \sigma(t \cdot g^{-1})(\xi) + d_t \sigma(t \cdot g^{-1}) \cdot d_x (g^{-1})(\xi) \right) \cdot (d_t g (t \cdot g^{-1}))^{-1} & +\\
& d_t g (t \cdot g^{-1}) \cdot \sigma(t \cdot g^{-1}) \cdot \left((d_x d_t g) (t \cdot g^{-1})(\xi) + (d_t d_t g)(t \cdot g^{-1}) \cdot d_x (g^{-1})(\xi)\right).
\end{array}
$$
The term $d_t g (t \cdot g^{-1}) \cdot d_x \sigma(t \cdot g^{-1}) \cdot (d_t g (t \cdot g^{-1}))^{-1}$ is equal to $g \ast (d_{dR} \sigma)$.

Secondly, we compute
$$
\begin{array}{lcl}
\omega'_{U}(-)(\xi)(u) & = & \tilde{\ast}'(\hat{\xi}_{u}(\psi')) = \tilde{\ast}'(({\psi'}^* \omega)(\xi_u)) = \tilde{\ast}'(\omega(d(\psi (t \cdot g(u)^{-1}))(\xi_u)))\\
& = & \tilde{\ast}'(\omega(d_x \psi (t \cdot g(u)^{-1}) + d_t (\psi (t \cdot g(u)^{-1}))(\xi_u)))\\
& = & \underbrace{\tilde{\ast}'(g(u)_* (\psi^* \omega(\xi_u)))}_{(a)} + \underbrace{\tilde{\ast}'(\omega(d_t (\psi (t \cdot g(u)^{-1}))(\xi_u)))}_{(b)},
\end{array}
$$
where the term $(a)$ is obtained by means of the definition of $\omega$ since
$$
\begin{array}{lcl}
\omega(d_x \psi (t \cdot g(u)^{-1}))(\xi_u) & = & -\alpha_{(u,\, \psi(t \cdot g(u)^{-1}))} \cdot s_{(u,\, \psi(t \cdot g(u)^{-1}))}^{-1} \cdot d_x \psi (t \cdot g(u)^{-1})(\xi_u)\\
& = & -\alpha_{(u,\, \psi(t \cdot g(u)^{-1}))} \cdot s_{(u,\, \psi(t))}^{-1} \cdot d_x \psi (t)(\xi_u)\\
& = & g(u)_*(-\alpha_{(u,\, \psi(t))} \cdot s_{(u,\, \psi(t))}^{-1} \cdot d_x \psi (t)(\xi_u))\\
& = & g(u)_*(\psi^* \omega(\xi_u)).
\end{array}
$$
We compute this term $(a)$ and the second term $(b)$ separately. The map $g_* -$ acts nontrivially only on the $W^0_m$ part of $W_m$. Therefore, we have the following computation: let $\psi_s$ be a path in $(\Rb^m)^{\coor}$ with tangent vector $s_{(0,\, \Id)} \cdot \psi^* \omega(\xi_u)$ at $s = 0$, we get that $g(u) \cdot \psi_s \cdot g(u)^{-1}$ has $s_{(0,\, \Id)} \cdot g(u)_* (\psi^* \omega(\xi_u))$ as a tangent vector at $s = 0$. The fact that $(g \cdot \psi_s \cdot g^{-1}) \ast (g \ast \sigma) = (g \cdot \psi_s \cdot g^{-1} \cdot g) \ast \sigma = (g \cdot \psi_s) \ast \sigma = g \ast (\psi_s \ast \sigma)$ shows that
$$
\tilde{\ast}'(g_* (\psi^* \omega(\xi)))(g \ast \sigma) = g \ast (\tilde{\ast}' (\psi^* \omega(\xi))(\sigma)).
$$
On the other hand, to calculate the term $(b)$, we remark that
$$
\begin{array}{lcl}
\omega(d_t (\psi (t \cdot g^{-1}))(\xi)) & = & -s_{(0,\, \Id)} \cdot d_t \psi(t \cdot g^{-1})^{\coor}(0,\, \Id)^{-1} \cdot d_t \psi (t \cdot g^{-1}) \cdot d_x (g^{-1})(\xi)\\
& = & -s_{(0,\, \Id)} \cdot d_g (g^{-1})^{-1} \cdot d_x (g^{-1})(\xi).
\end{array}
$$
Therefore $\tilde{\ast}'(\omega(d_t (\psi (t \cdot g^{-1}))(\xi)))(g \ast \sigma(t)) = - \ast' (d_g (g^{-1})^{-1} \cdot d_x (g^{-1})(\xi))(g \ast \sigma(t))$. Moreover, the vector field $d_g (g^{-1})^{-1} \cdot d_x (g^{-1})(\xi)$ has no constant constant part with respect to the variable $t$ so we can compute
$$
(d_{t} A)(t \cdot g^{-1}) \cdot d_g (g^{-1}) (d_g (g^{-1})^{-1} \cdot d_x (g^{-1})(\xi)) = (d_{t} A)(t \cdot g^{-1}) \cdot d_x (g^{-1})(\xi),
$$
and
$$
d_g (d_t g)(d_g (g^{-1})^{-1} \cdot d_x (g^{-1})(\xi)) = d_t (d_x g (\xi)) = d_x d_t g (\xi).
$$
It follows from the remark coming after Definition \ref{rem1} that\\
$\tilde{\ast}'(\omega(d_t (\psi (t \cdot g^{-1}))(\xi)))(g \ast \sigma(t)) = $
$$
\begin{array}{clc}
- & \left((d_x d_t g) (t \cdot g^{-1})(\xi) + (d_t d_t g)(t \cdot g^{-1}) \cdot d_x (g^{-1})(\xi)\right) \cdot \sigma(t \cdot g^{-1}) \cdot (d_t g (t \cdot g^{-1}))^{-1}\\
- & d_t g (t \cdot g^{-1}) \cdot d_t \sigma(t \cdot g^{-1}) \cdot d_x (g^{-1})(\xi) \cdot (d_t g (t \cdot g^{-1}))^{-1}\\
- & d_t g (t \cdot g^{-1}) \cdot \sigma(t \cdot g^{-1}) \cdot \left((d_x d_t g) (t \cdot g^{-1})(\xi) + (d_t d_t g)(t \cdot g^{-1}) \cdot d_x (g^{-1})(\xi)\right)
\end{array}
$$
and the Lemma is proved.
$\cqfd$
\end{pf}

\begin{dei}
We consider now on $\tE$ the left-action of $\GL_m(\Rb) \subset G^0_m$ given by
$$
(g,\, (x,\, \varphi,\, A(t))) \mapsto (\rho(g,\, (x,\, \varphi)),\, g \ast A(t)) = (x,\, \varphi \cdot g^{-1},\, g \ast A(t)),
$$
where $g$ is seen as a linear automorphism of $\R^m_{\mathrm{for}}$. We define by $\bE$ the quotient space
$$
\bE := \Mco \times_{\GL_m(\Rb)} \Mrm_m \left( \llbracket t^1,\, \ldots ,\, t^m \rrbracket \right).
$$
It is a fiber bundle over $\Maf := \Mco / \GL_m (\Rb)$.
\end{dei}
The fibers of the fiber bundle $\Maf$ over $M$ are contractible. Because of this, there exists at least one section $M \rightarrow \Maf$. We fix such a section $\psi^{\mathrm{aff}}$ and we denote the pullback bundle on $M$ by
$$
E := (\psi^{\mathrm{aff}})^* \bE \rightarrow M.
$$
The restriction of $\psi^{\mathrm{aff}}$ on a contractible open set $U$ of $M$ is an equivalence class of sections $\psi : U \rightarrow U^{\coor}$ and two sections $\psi$ and $\psi'$ are equivalent if there exists a smooth map $g : U \rightarrow \GL_m(\Rb)$ such that $\psi'(u)(t) = \psi(u) (t \cdot g(u)^{-1})$. Restricted to $U$, a section $\sigma$ of $E$ is an equivalence class of sections $\sigma_U : U \rightarrow E_{U} = U \times \mathrm{M}_{m}\left(\Rb \llbracket t^1,\, \ldots ,\, t^m \rrbracket \right)$, each being associated to a section $U \rightarrow U^{\coor}$. Two such sections $\sigma_U$ and $\sigma'_U$, associated to $\psi$ and $\psi'$ respectively, are equivalent if there exists a smooth map $g : U \rightarrow \GL_m(\Rb)$ such that $\psi'(u)(t) = \psi(u) (t \cdot g(u)^{-1})$ and $\sigma'_U = g \ast \sigma_U$. Let $\sigma$ associated to $\psi$. We define $(\nabla \sigma)_{|U} := \nabla_U(\psi)(\sigma_U)$.

\begin{thm}\label{connglob}
The map $\nabla : \Gamma(E) \rightarrow \Omega_{M}^{1} \otimes_{\C^{\infty}(M)} \Gamma(E)$ is a well-defined connection on $E$.
\end{thm}

\begin{pf}
Lemma \ref{comp} shows that the connection $\nabla$ is well-defined on equivalence classes of sections of $E_U$ for any contractible open set $U$. It therefore induces a globally defined connection on $E$.
$\cqfd$
\end{pf}

\subsection{Complex manifold and flat section}\label{cplxflat}

We assume in this section that a section $\psi^{\mathrm{aff}} : M \rightarrow \Maf$ is fixed and we denote by $\nabla$ the associated connection on $E = (\psi^{\mathrm{aff}})^* \bE$. A complex structure on a formal manifold of dimension $m$ is given by an element in $\mathrm{M}_{m}\left(\Rb \llbracket t^1,\, \ldots ,\, t^m \rrbracket \right)$ and we can therefore see the set of $\Cx_{\infty}$-algebras on $\Rb^{m}$ as a subset of $\mathrm{M}_{m}\left(\Rb \llbracket t^1,\, \ldots ,\, t^m \rrbracket \right)$. Let $J$ be a smooth complex structure on the manifold $M$. Associated to any local coordinates system $\tilde{\varphi}$ around a point $x = \tilde{\varphi}(0)$, the Taylor series of $\tilde{\varphi}_{*}(J) := (d \tilde{\varphi})^{-1} \cdot J(\tilde{\varphi}) \cdot d \tilde{\varphi}$ at $0$ provides a matrix $A(J,\, \tilde{\varphi}) \in \{ \mathrm{Cx}_{\infty}\textrm{-algebras}\}$, which depends only on the restriction $\varphi$ of $\tilde{\varphi}$ to $\R^{m}_{\mathrm{for}}$. We therefore denote it by $A(J,\, \varphi)$. We have the relation $A(J,\, \varphi \cdot g^{-1}) = g \ast A(J,\, \varphi)$ for any automorphism $g$ of $\R^{m}_{\mathrm{for}}$, so it follows that a complex structure on $M$ provides a section of the fiber bundle $E_{cx}$ over $M$
$$
E_{cx} := (\psi^{\mathrm{aff}})^* (M^{\textrm{coor}} \times_{\GL_m(\Rb)} \{ \mathrm{Cx}_{\infty}\textrm{-algebras}\}) \subseteq E.
$$
In the sequel, we characterize sections of $E_{cx}$ which come from smooth complex structures on $M$.\\

Let $\tilde{\sigma} : TM \rightarrow TM$ be a smooth function. As just seen, we can associate to any such $\tilde{\sigma}$ a smooth section $\sigma : M \rightarrow E$ of $E$. We call \emph{geometric smooth functions} the sections of $E$ defined in this manner. We are interested in the connection $\nabla$ because of the following result.

\begin{prop}\label{flat}
Flat sections for the connection $\nabla$, that is, sections $\sigma$ such that $\nabla \sigma = 0$, are geometric smooth functions $M \rightarrow E$.
\end{prop}

\begin{pf}
Restricted to some contractible open subset $U$ of $M$, the section $\psi^{\mathrm{aff}}$ is an equivalence class of section $U \rightarrow U^{\coor}$. We fix a section $\psi$ in this class. A section $\sigma$ of $E$ is a class of sections of $E_U$. We denote by $\sigma_U : U \rightarrow E_{U} = U \times \mathrm{M}_{m}\left(\Rb \llbracket t^1,\, \ldots ,\, t^m \rrbracket \right)$ the representative section associated to $\psi$. Because of the definition of $\nabla$, the section $\sigma$ is flat for the connection $\nabla$ if and only if the section $\sigma_U$ is flat for the connection $\nabla_U$.

Let $u_0$ be a  point in $U$ and $\tilde{\psi}_{u_0} : \Rb^m \rightarrow U$ a diffeomorphism whose restriction to $\R^m_{\mathrm{for}}$ is $\psi_{u_0}$ (the target space might be smaller than $U$ and we replace $U$ by this smaller contractible open set in this case). We have a sequence of maps
$$
\R^m_{\mathrm{for}} \xrightarrow{\psi_u} (U,\, \C^{\infty}_U) \xrightarrow{\tilde{\psi}_{u_0}^{-1}} \R^m \xleftarrow{\lambda_{\tilde{\psi}_{u_0}^{-1}(u)}} \R^m_{\mathrm{for}},
$$
where the last map, which sends $0$ to $\tilde{\psi}_{u_0}^{-1}(u)$ and a function on $\Rb^m$ to its Taylor series at the point $\tilde{\psi}_{u_0}^{-1}(u)$, is an immersion of locally ringed spaces. It is therefore invertible on its image and we obtain for all $u \in U$ an automorphism $g(u) : \R^m_{\mathrm{for}} \rightarrow \R^m_{\mathrm{for}}$. The map $g : U \rightarrow G^0_m$ is smooth and by Lemma \ref{comp}, we have that $\nabla_U(\psi) (\sigma_U) = 0$ if and only if $\nabla_U(\psi') (g \ast \sigma) = 0$ for $\psi'(u) = \tilde{\psi}_{u_0} \cdot \lambda_{\tilde{\psi}_{u_0}^{-1}(u)}$. We now compute $\nabla_U(\psi') (g \ast \sigma) (\xi)$ using the chart given by $\tilde{\psi}_{u_0}$. We get
$$
d_{dR} (g \ast \sigma)(\xi) = d_{dR} \left( (g \ast \sigma) \cdot {\tilde{\psi}_{u_0}} \right) \left((d {\tilde{\psi}_{u_0}})^{-1} (\xi) \right)
$$
and
$$
\begin{array}{lcl}
\omega'_{U} (-)(\xi) & = & \tilde{\ast}' (\hat{\xi}(\psi')) =  \tilde{\ast}' (\omega(d_t \tilde{\psi}_{u_0} \cdot d_u \lambda_{\tilde{\psi}_{u_0}^{-1}(u)} (\xi)))\\
& = & \tilde{\ast}' \left(-s_{(0,\, \Id)}^{-1} \cdot d_t \left(\tilde{\psi}_{u_0} \cdot \lambda_{\tilde{\psi}_{u_0}^{-1}(u)} \right)^{-1} \cdot (d_t \tilde{\psi}_{u_0} \cdot d_u \lambda_{\tilde{\psi}_{u_0}^{-1}(u)} (\xi)) \right)\\
& = & - \tilde{\ast}' \left(s_{(0,\, \Id)}^{-1} \cdot d_u \lambda_{\tilde{\psi}_{u_0}^{-1}(u)} (\xi)) \right) = - \tilde{\ast}' \left(d_t \tilde{\psi}_{u_0}^{-1} (\xi) \right).
\end{array}
$$
The vector $d_t \tilde{\psi}_{u_0}^{-1} (\xi_{u_0})$ is a constant vector field on $\R^m_{\mathrm{for}}$, hence it acts through $\tilde{\ast}'$ by differentiation on $g \ast \sigma$. When $\xi$ varies, the vectors $d_t \tilde{\psi}_{u_0}^{-1} (\xi_{u_0})$ cover $\Rb^m$ and we obtain that the equality $\nabla_U(\psi') (g \ast \sigma) = 0$ is equivalent to
\begin{equation}\label{eqsmooth}
\frac{1}{\underline{\alpha}!} \frac{\partial}{\partial u^{\alpha}} \left(\sigma_{\alpha_1,\, \ldots,\, \alpha_n,\, a}^b (u)\right) = \frac{1}{\underline{\alpha'}!} \sigma_{\alpha_1,\, \ldots ,\, \alpha_n,\, \alpha,\, a}^b (u), \textrm{ for all } \alpha,\, \alpha_1,\, \ldots ,\, \alpha_n,\, a,\, b,
\end{equation}
where $(g \ast \sigma)(u) = \sum_{n \geq 0} \sigma_{\alpha_1,\, \ldots ,\, \alpha_n,\, a}^b(u) \frac{1}{\underline{\alpha}!} t^{\alpha_1} \cdots t^{\alpha_n} \gamma^a \partial_b$ and $\underline{\alpha'} := \{ \alpha_1,\, \ldots ,\, \alpha_k,\, \alpha \}$. This is equivalent to the fact that the matrix $(g \ast \sigma) (u)(t) \in \Mrm_m \left( \Rb \llbracket t^1,\, \ldots ,\, t^m \rrbracket \right)$ is the infinite jet of the smooth function $(g \ast \sigma) (-)(0) = \left(\sigma_a^b (-) \right)_{a,\, b} : U \rightarrow \Mrm_m (\Rb)$. We finally obtain that the section $\sigma$ is a geometric smooth function if and only if it is a flat section for the connection $\nabla$.
$\cqfd$
\end{pf}

As a corollary, we get the following theorem

\begin{thm}\label{cxflatthm}
Complex structures on a smooth manifold $M$ correspond to flat sections of the fiber bundle $E_{cx} = (\psi^{\mathrm{aff}})^* (M^{\textrm{coor}} \times_{\GL_m(\Rb)} \{ \mathrm{Cx}_{\infty}\textrm{-algebras}\})$.
\end{thm}

\begin{pf}
We have already seen at the beginning of this section that a complex structure on $M$ provides a section of $E_{cx}$. By Proposition \ref{flat}, the associated section is flat. Conversely, by Proposition \ref{flat}, a flat section of $E_{cx}$ is a geometric smooth function hence corresponds to a smooth endomorphism $J$ of $TM$. The property to be a complex structure (almost complex structure and integrability condition) can be read pointwise by means of the Taylor coefficient of $J$ of order $0$ and $1$. Because of the fact that we have considered a section of $E_{cx}$ and by Theorem \ref{equiv1}, we get that $J$ is a smooth complex structure.
$\cqfd$
\end{pf}

\subsection{Holomorphic map and flat section}\label{holoflat}

Sections \ref{coorandconn} and \ref{cplxflat} extend to the case of holomorphic maps between complex manifolds in the following manner. Let $(M,\, J_{M})$ and $(N,\, J_{N})$ be two complex manifolds of dimension $m$ and $n$ respectively. We denote by $(M \times N)^{\coor}$ the manifold
$$
(M \times N)^{\textrm{coor}} := \left\{ \begin{array}{c}
(x,\, y,\, \varphi,\, \psi) \textrm{ such that } x\in M,\, y \in N \textrm{ and}\\
\varphi : \mathcal{R}^{m}_{\mathrm{for}} \rightarrow (\mathcal{M},\, x) \textrm{ and } \psi : \mathcal{R}^{n}_{\mathrm{for}} \rightarrow (\mathcal{N},\, y)\\
\textrm{are pointed immersions of loc. ringed spaces} \end{array} \right\},
$$
where $\M = (M,\, \C^{\infty}_{M})$ and $\N = (N,\, \C^{\infty}_{N})$. It is a fiber bundle over $M \times N$.
\begin{dei}
We define the trivial vector bundle (over $(M \times N)^{\coor}$)
$$
\tilde{F} := (M \times N)^{\coor} \times \Map(\R^m_{\mathrm{for}},\, \R^n_{\mathrm{for}}),
$$
where $\Map(\R^m_{\mathrm{for}},\, \R^n_{\mathrm{for}}) \cong \left(\Rb \llbracket t^1,\, \ldots ,\, t^m \rrbracket^{\geq 1} \right)^n$. The group $G^0_m \times G^0_n$ acts on $\tilde{F}$ on the left by\\
$((g,\, h),\, (x,\, y,\, \varphi,\, \psi,\, v(t))) \xmapsto{(\delta,\, \diamond)}$
$$
(\delta((g,\, h),\, (x,\, y,\, \varphi,\, \psi)),\, (g,\, h) \diamond v(t)) := (x,\, y,\, \varphi \cdot g^{-1},\, \psi \cdot h^{-1},\, h \cdot v(t \cdot g^{-1})).
$$
\end{dei}
We can derive the action $\diamond$ and extend it to get a map
$$
\tilde{\diamond}' : W_m \times W_n \rightarrow \chi \left(N \times \left(\Rb \llbracket t^1,\, \ldots ,\, t^m \rrbracket^{\geq 1} \right)^n \right),
$$
The differentiation of the action $\diamond$ gives a map $W^0_m \times W^0_n$ to $\chi \left( \left(\Rb \llbracket t^1,\, \ldots ,\, t^m \rrbracket^{\geq 1} \right)^n \right)$, the constant vector fields $\zeta_m$ in $W_m$ act by differentiation of the vector in $\left(\Rb \llbracket t^1,\, \ldots ,\, t^m \rrbracket^{\geq 1} \right)^n$ with respect to the variables $t^k$, that is, $(y,\, v(t))$ is sent to $(0,\, d_t v(\zeta_m))$, and the constant vector fields $\zeta_n$ in $W_n$ send $(y,\, v(t))$ to $(v(0)(\zeta_n),\, 0)$.\\

Let $U$ be an open subset of $M$ and $V$ be an open subset of $N$. Associated to a section $\Psi : U \rightarrow (U \times V)^{\coor}$, we define a connection on the fiber bundle $F_U := U \times V \times \left(\Rb \llbracket t^1,\, \ldots ,\, t^m \rrbracket^{\geq 1} \right)^n$ over $U$ in the same way as in Proposition \ref{conn}. Lemma \ref{comp} extends to this setting. The subgroup $\GL_m (\Rb) \times \GL_n (\Rb)$ of $G^0_m \times G^0_n$ acts on $\tilde{F}$ and we define the quotient space $\overline{F} := (M \times N)^{\coor} \times_{\GL_m (\Rb) \times \GL_n (\Rb)} \left(\Rb \llbracket t^1,\, \ldots ,\, t^m \rrbracket^{\geq 1} \right)^n$. It is a fiber bundle over $(M \times N)^{\textrm{aff}} := (M \times N)^{\coor} / \left( \GL_m (\Rb) \times \GL_n (\Rb) \right)$. Moreover, $(M \times N)^{\textrm{aff}}$ is a fiber bundle over $M \times N$ whose fibers are contractible. We fix a section $\Psi^{\textrm{aff}} : M \times N \rightarrow (M \times N)^{\mathrm{aff}}$ and we denote by $F := (\Psi^{\mathrm{aff}})^* \overline{F}$ the pullback bundle on $M \times N$, that we see as a fiber bundle on $M$. Theorem \ref{connglob} and Proposition \ref{flat} extend to this setting. 
\begin{dei}
An $\infty$-morphism between $\Cx_{\infty}$-algebras $\Rb^m$ and $\Rb^n$ is an element in
$$
\left(\Rb \llbracket t^1,\, \ldots ,\, t^m \rrbracket^{\geq 1} \right)^n,
$$
thus we can define the fiber bundle
$$
F_{cx} := (\Psi^{\mathrm{aff}})^* \left( (M \times N)^{\coor} \times_{\GL_m (\Rb) \times \GL_n (\Rb)} \{ \infty \textrm{-morphisms}\} \right)
$$
over $M$.
\end{dei}
The analog of Theorem \ref{cxflatthm} holds and we obtain eventually the following theorem.

\begin{thm}\label{equiv2}
Let $M$ and $N$ be two smooth manifolds. We write $E_{cx}(M)$ and $F_{cx}(M,\, N)$ to emphasize the fact that these fiber bundles, previously defined, depend on $M$ and $N$. There is an equivalence of categories
$$\begin{array}{ccc}
\left\{ \begin{gathered} \textrm{Flat sections of } E_{cx}(M) \textrm{ with} \\ \textrm{flat sections of } F_{cx}(M,\, N) \end{gathered} \right\} & \xrightarrow{\cong} & \left\{ \begin{gathered} \textrm{Complex structures on } M \textrm{ with}\\ \textrm{holomorphic maps between } M \textrm{ and } N \end{gathered} \right\}.
\end{array}$$
\end{thm}\

\begin{appendix}

\section{Operadic decomposition maps}\label{Lie}

The decomposition maps for the operads $\As^{\ash}$ and $\Lie^{\ash}$ are described in Sections 9.1.5 and 10.1.6 in \cite{LodayVallette}. The formulas for $\As_{1}^{\ash}$ and $\Lie_{1}^{\ash}$ are similar, the only differences are degrees and signs.

\subsection{Associative case}

The operad encoding associative algebras endowed with a product of cohomological degree 1 is given by
$$
\As_{1} := \left. \T(s^{-1}E_{A}) \middle/ \left(s^{-2}R_{A} \right) \right. ,
$$
where $E_{A}$ is the free $\Sb$-module generated in arity $2$ by an element $\overline{\mu}$, or $\as$, and $R_{A}$ is the free $\Sb$-module generated in arity $3$ by the associativity relation $(\overline{\mu}\ ;\ \overline{\mu},\, \Id) + (\overline{\mu}\ ;\ \Id,\, \overline{\mu})$ or $\ass$.

\begin{rem}
The ``$+$'' sign in the associativity relation gives that $A$ is a $\As_{1}$-algebra if, and only if, $sA$ is an associative algebra (in the classical sense).
\end{rem}

Similarly to the case of the operad $\As$ encoding associative algebras, we can compute the Koszul dual cooperad to the operad $\As_1$.

\begin{prop}
The Koszul dual cooperad \emph{$\As^{\ash}_{1}$} is generated (as an $\Sb$-module) in arity $n$ by the element
$$
\overline{\mu}_{n}^{c} := \sum_{t\, \in PBT_{n}} t,
$$
where $PBT_{n}$ is the set of planar binary trees on $\overline{\mu}$ with $n$ leaves. It follows that the infinitesimal decomposition map on \emph{$\As^{\ash}_{1}$} is given by
\emph{$$
\Delta^{(1)}_{\As_{1}^{\ash}}(\overline{\mu}_{n}^{c}) = \sum_{\substack{l+q+r=n+1\\ p=l+r+1\geq 1,\, q \geq 1}} (\overline{\mu}_{p}^{c}\, ;\, \underbrace{\Id,\, \ldots ,\, \Id}_{l} ,\, \overline{\mu}_{q}^{c},\, \underbrace{\Id ,\, \ldots ,\, \Id}_{r}),
$$}
and that the full decomposition map is given by
\emph{$$
\Delta_{\As_{1}^{\ash}}(\overline{\mu}_{n}^{c}) = \sum_{q_{1}+\cdots +q_{p}= n} (\overline{\mu}_{p}^{c}\, ;\, \overline{\mu}_{q_{1}}^{c},\, \ldots ,\, \overline{\mu}_{q_{p}}^{c}).
$$}
\end{prop}

\subsection{Lie case}

The operad encoding Lie algebras endowed with a bracket of cohomological degree 1 is given by
$$
\Lie_{1} := \left. \T(s^{-1}E_{L}) \middle/ \left(s^{-2}R_{L} \right) \right. ,
$$
where $E_{L}$ is the $\Sb$-module generated in arity $2$ by a symmetric element $\lie$ and $R_L$ is the $\Sb$-module generated by the Jacobi relation: $\jac$.\\

There is a morphism of operads $\Lie_{1} \rightarrow \As_{1}$ defined by
$$
\slie \mapsto s^{-1}\overline{\mu} + s^{-1}\overline{\mu}^{(12)}.
$$
It is well-defined since it sends the Jacobi relation to a linear combination of associativity relations. For the same reason, there is a morphism of cooperads $\Lie_{1}^{\ash} \rightarrow \As_{1}^{\ash}$ defined by
$$
\lie \mapsto \overline{\mu} + \overline{\mu}^{(12)}.
$$

We make use of the associative case to compute the Koszul dual cooperad associated to $\Lie_1$.

\begin{prop}
The Koszul dual cooperad $\Lie_{1}^{\ash}$ is 1-dimensional in each arity, \emph{$\Lie_{1}^{\ash}(n) \cong \Kb \cdot \overline{l}_{n}^{c}$}, where $\overline{l}_{n}^{c}$ is an element of degree $0$ such that
$$
\overline{l}_{n}^{c} := \sum_{\sigma \in \Sb_{n}} (\overline{\mu}_{n}^{c})^{\sigma}.
$$
(It follows that $\Sb_n$ acts trivially on $\overline{l}_n^c$.) The infinitesimal decomposition map on \emph{$\Lie_{1}^{\ash}$} is given by
\emph{$$
\Delta^{(1)}_{\Lie_{1}^{\ash}}(\overline{l}_{n}^{c}) = \sum_{\substack{p+q=n+1\\ p,\, q \geq 1}} \sum_{\sigma \in Sh_{q,\, p-1}^{-1}} (\overline{l}_{p}^{c} \circ_{1} \overline{l}_{q}^{c})^{\sigma},
$$}
where $Sh_{q,\, p-1}^{-1}$ is the set of $(q,\, p-1)$-unshuffles, that is, inverses of $(q,\, p-1)$-shuffles. The formula for the full decomposition map is given by
\emph{$$
\Delta_{\Lie_{1}^{\ash}}(\overline{l}_{n}^{c}) = \sum_{\substack{\{q_{1},\, \ldots ,\, q_{p} \}\\ q_{1}+\cdots +q_{p}= n}} \sum_{\sigma \in Sh_{q_{1},\, \ldots ,\, q_{p}}^{-1}} \frac{1}{N_{\underline{q}}} (\overline{l}_{p}^{c}\, ;\, \overline{l}_{q_{1}}^{c},\, \ldots ,\, \overline{l}_{q_{p}}^{c})^{\sigma},
$$}
where $N_{\underline{q}} := \prod_{i=1}^{m} \max(n_{i},\, 1)$ with $\underline{q} = \{ q_{1},\, \ldots ,\, q_{n}\} = \{ \underbrace{1,\, \ldots ,\, 1}_{n_{1} \mathrm{\ times}},\, 2,\, \ldots ,\, \underbrace{m,\, \ldots ,\, m}_{n_{m} \mathrm{\ times}} \}$.
\end{prop}

\section{Distributive laws and decomposition map}\label{distri}

In this appendix, we define distributive laws for cooperads in order to compute the decomposition map of the Koszul dual cooperad of an operad endowed with a distributive law. We dualize the presentation given by Loday and Vallette \cite{LodayVallette}, Section 8.6. We emphasize however that we work here with cohomological degree and not with homological degree. We will always consider the opposite of the signs appearing in \cite{LodayVallette} and the chain complexes will be bounded above.

\subsection{Distributive law for cooperads}\label{distrilaw}

Let $(\C,\, \Delta_{\C},\, \epsilon_{\C})$ and $(\D,\, \Delta_{\D},\, \epsilon_{\D})$ be two cooperads. A morphism of $\Sb$-modules $\Lambda^{c} : \D \circ \C \rightarrow \C \circ \D$ is called a \emph{distributive law for cooperads} if the following diagrams are commutative:
$$
\mathrm{(I)} \hspace{1cm} \xymatrix{\D \circ \C \ar[rr]^{\Lambda^{c}} \ar[d]_{\Delta_{\D} \circ \Id_{\C}} & & \C \circ \D \ar[d]^{\Id_{\C} \circ \Delta_{\D}}\\
\D \circ \D \circ \C \ar[r]^{\Id_{\D} \circ \Lambda^{c}} & \D \circ \C \circ \D \ar[r]^{\Lambda^{c} \circ \Id_{\D}} & \C \circ \D \circ \D,}
$$
$$
\mathrm{(II)} \hspace{1cm} \xymatrix{\D \circ \C \ar[rr]^{\Lambda^{c}} \ar[d]_{\Id_{\D} \circ \Delta_{\C}} & & \C \circ \D \ar[d]^{\Delta_{\C} \circ \Id_{\D}}\\
\D \circ \C \circ \C \ar[r]^{\Lambda^{c} \circ \Id_{\C}} & \C \circ \D \circ \C \ar[r]^{\Id_{C} \circ \Lambda^{c}} & \C \circ \C \circ \D,}
$$
$$
\mathrm{(i)} \hspace{0.5cm} \xymatrix{& \C &\\
\D \circ \C \ar[ur]^{\epsilon_{\D} \circ \Id_{\C}} \ar[rr]^{\Lambda^{c}} & & \C \circ \D \ar[ul]_{\Id_{\C} \circ \epsilon_{\D}},}
\hspace{1cm} \mathrm{(ii)} \hspace{0.5cm} \xymatrix{& \D &\\
\D \circ \C \ar[ur]^{\Id_{\D} \circ \epsilon_{\C}} \ar[rr]^{\Lambda^{c}} & & \C \circ \D \ar[ul]_{\epsilon_{\C} \circ \Id_{\D}}.}
$$

\begin{propA}
If $\Lambda^{c} : \D \circ \C \rightarrow \C \circ \D$ is a distributive law for the cooperads $\C$ and $\D$, then $\D \circ \C$ is a cooperad for the decomposition map
$$
\Delta_{\Lambda^{c}} := (\Id_{\D} \circ \Lambda^{c} \circ \Id_{\C}) (\Delta_{\D} \circ \Delta_{\D}) : \D \circ \C \rightarrow (\D \circ \C) \circ (\D \circ \C),
$$
and for the counit
$$
\epsilon_{\Lambda^{c}} := \epsilon_{\D} \circ \epsilon_{\C} : \D \circ \C \rightarrow \Im.
$$
\end{propA}

\begin{pf}
It is enough to dualize the proof of Proposition 8.6.2 in \cite{LodayVallette}. To simplify the notations, we write $\C \D$ instead of $\C \circ \D$. The following diagram commutes
$$
\xymatrix@C=8pt@R=16pt{\D \C \ar[rrrr] \ar[dd] &&&& \D \D \C \C \ar[d] \ar[rr] && \D \C \D \C \ar[d]\\
&& \textrm{ coassoc. of } \Delta_{\D},\, \Delta_{\C} && \D \D \D \C \C \ar[d] \ar[rr] && \D \D \C \D \C \ar[dd]\\
\D \D \C \C \ar[rr] \ar[dd] && \D \D \C \C \C \ar[rr] && \D \D \D \C \C \C \ar[d] & \mathrm{(II)} &\\
&& \mathrm{(I)} && \D \D \C \D \C \C \ar[rr] \ar[d] && \D \D \C \C \D \C \ar[d] \\
\D \C \D \C \ar[rr] && \D \C \D \C \C \ar[rr] && \D \C \D \D \C \C \ar[rr] && \D \C \D \C \D \C \D \C.}
$$
The arrows in this diagram are composite product of identities, $\Lambda^{c}$, $\Delta_{\D}$ and $\Delta_{\C}$, for example, $\Delta_{\D} \circ \Delta_{\C} : \D \C \rightarrow \D \D \C \C$, and are uniquely determined by their source and their target. (Remember that $\Delta_{\C}$ and $\Delta_{\D}$ are coassociative.) The two empty squares commute because the composite product $\circ$ is a bifunctor. The counit property is proved in a similar way by means of (i) and (ii).
$\cqfd$
\end{pf}

\subsection{Decomposition map by means of distributive law}\label{naturali}

Let $(\A,\, \gamma_{\A},\, \iota_{\A})$ and $(\B,\, \gamma_{\B},\, \iota_{\B})$ be two operads. We assume that $\A$ and $\B$ are quadratic with quadratic presentations $\A = \P(V,\, R) := \T (V)/(R)$ and $\B = \P(W,\, S) := \T (W)/(S)$ and that we have a rewriting rule $\lambda : W \circ_{(1)} V \rightarrow V \circ_{(1)} W$. We denote the graph of $\lambda$ by
$$
D_{\lambda} := \langle T - \lambda (T),\, T \in W \circ_{(1)}  V \rangle \subset \T (V \oplus W)^{(2)}.
$$
Let $\A^{\ash} = \C (sV,\, s^{2}R)$, resp. $\B^{\ash} = \C (sW,\, s^{2}S)$, be the Koszul dual cooperad of $\A$, resp. $\B$. The categorical coproduct of the two operads, denoted $\A \vee \B$, is equal to $\T (V \oplus W)/(R\oplus S)$. On the other side, the categorical coproduct of the two Koszul dual cooperads, denoted $\A^{\ash} \vee \B^{\ash}$, is equal to $\C (sV\oplus sW,\, s^{2}R \oplus s^{2}S)$. We remark therefore that $(\A \vee \B)^{\ash} = \A^{\ash} \vee \B^{\ash}$.

We now consider the cooperad $\A^{\ash} \vee^{\lambda} \B^{\ash}$ given by the following presentation
$$
\A^{\ash} \vee^{\lambda} \B^{\ash} := \C (sV\oplus sW,\, s^{2}R \oplus s^{2}D_{\lambda} \oplus s^{2}S).
$$
The categorical coproduct of $\A^{\ash}$ and $\B^{\ash}$ injects itself in this cooperad. Moreover, we can build the following map
$$
\C (sV \oplus sW,\, s^{2}R \oplus s^{2}D_{\lambda} \oplus s^{2}S) \hookrightarrow \T^{c}(sV \oplus sW) \twoheadrightarrow \T^{c}(sW) \circ \T^{c}(sV),
$$
where the second arrow is the projection $p_{1}$ which sends any tree in $sV$ and $sW$ containing a subtree in $sV \circ_{(1)} sW$ to $0$ and is identity on other trees. This composition factors through the inclusion $\C (sV,\, s^{2}R) \circ \C (sW,\, s^{2}S) \hookrightarrow \T^{c}(sV) \circ \T^{c}(sW)$ to give a morphism of $\Sb$-modules
$$
i_{1} : \A^{\ash} \vee^{\lambda} \B^{\ash} \rightarrow \B^{\ash} \circ \A^{\ash}.
$$
The map $i_{1}$ is an inclusion since the composition $D_{\lambda} \hookrightarrow (W \circ_{(1)} V) \oplus (V \circ_{(1)} W) \twoheadrightarrow W \circ_{(1)} V$ is an inclusion. Indeed, given an element in $\B^{\ash} \circ \A^{\ash}$, there is at most one possibility to build an element in $\A^{\ash} \vee^{\lambda} \B^{\ash}$ from it. (The composition is also a surjection however $i_{1}$ is not necessarily a surjection.) Similarly, we get a morphism of $\Sb$-modules $i_{2} : \A^{\ash} \vee^{\lambda} \B^{\ash} \rightarrow \A^{\ash} \circ \B^{\ash}$, which is neither necessarily an inclusion, nor necessarily a surjection.

To make this morphism $i_{1}$ easier to understand, we describe it partially. We recall that
$$
\A^{\ash} \vee^{\lambda} \B^{\ash} = \C(sV \oplus sW,\, s^{2}R \oplus s^{2}D_{\lambda} \oplus s^{2}S) = \Im \oplus sV \oplus sW \oplus s^{2}R \oplus s^{2}D_{\lambda} \oplus s^{2}S \oplus \cdots.
$$
Applying $i_{1}$, $\Im$ is sent on $\Im$ by the identity map and $sV$, resp. $sW$, is sent to $\Im \circ sV$, resp. $sW \circ I$, and $s^{2}R$, resp. $s^{2}S$, is sent to $\Im \circ s^{2}R$, resp. $s^{2}S \circ I$, and $s^{2}D_{\lambda}$ is sent to $i_{1}(s^{2} D_{\lambda}) \subseteq sW \circ_{(1)} sV \subset \B^{\ash} \circ \A^{\ash}$. The application $i_{2}$ is defined similarly.

\begin{rem}
The map $p : \A \circ \B \rightarrow \A \vee_{\lambda} \B$ given in \cite{LodayVallette}, Section 8.6.2, is defined dually.
\end{rem}

\begin{prop}\label{decomposition map}
Let $\A = \P(V,\, R)$ and $\B = \P(W,\, S)$ be two quadratic operads and $\A^{\ash} = \C(sV,\, s^{2}R)$ and $\B^{\ash} = \C(sW,\, s^{2}S)$ be their Koszul dual cooperads. For any morphism of $\Sb$-modules $\lambda : W \circ_{(1)} V \rightarrow V \circ_{(1)} W$ such that $i_{1} : \A^{\ash} \vee^{\lambda} \B^{\ash} \hookrightarrow \B^{\ash} \circ \A^{\ash}$ is an isomorphism, the composite
$$
\Lambda^{c} : \B^{\ash} \circ \A^{\ash} \xrightarrow{i_{1}^{-1}} \A^{\ash} \vee^{\lambda} \B^{\ash} \xrightarrow{i_{2}} \A^{\ash} \circ \B^{\ash}
$$
induces a distributive law for cooperads and a decomposition map
$$
\Delta_{\Lambda^{c}} := (\Id_{\B^{\ash}} \circ \Lambda^{c} \circ \Id_{\A^{\ash}}) \cdot (\Delta_{\B^{\ash}} \circ \Delta_{\A^{\ash}}) : \B^{\ash} \circ \A^{\ash} \rightarrow (\B^{\ash} \circ \A^{\ash}) \circ (\B^{\ash} \circ \A^{\ash}),
$$
and a counit
$$
\epsilon_{\Lambda^{c}} := \epsilon_{\B^{\ash}} \circ \epsilon_{\A^{\ash}} : \B^{\ash} \circ \A^{\ash} \rightarrow \Im,
$$
such that the map $i_{1} : \A^{\ash} \vee^{\lambda} \B^{\ash} \rightarrow (\B^{\ash} \circ \A^{\ash},\, \Delta_{\Lambda^{c}},\, \epsilon_{\Lambda^{c}})$ is an isomorphism of cooperads.
\end{prop}

\begin{pf}
The proof is dual to the proof of Proposition 8.6.4 in \cite{LodayVallette}. The previous descriptions of the application $i_{1}$ and $i_{2}$ provide the commutative diagrams
$$
\xymatrix@R=8pt{\B^{\ash} \circ \A^{\ash} \ar[rd]_{\epsilon_{\B^{\ash}} \circ \Id_{\A^{\ash}}} & \A^{\ash} \vee^{\lambda} \B^{\ash} \ar[l]_{i_{1}} \ar[r]^{i_{2}} & \A^{\ash} \circ \B^{\ash} \ar[dl]^{\Id_{\A^{\ash}} \circ \epsilon_{\B^{\ash}}}\\ & \A^{\ash} &} \textrm{ and } \xymatrix@R=8pt{\A^{\ash} \circ \B^{\ash} \ar[rd]_{\epsilon_{\A^{\ash}} \circ \Id_{\B^{\ash}}} & \A^{\ash} \vee^{\lambda} \B^{\ash} \ar[l]_{i_{2}} \ar[r]^{i_{1}} & \B^{\ash} \circ \A^{\ash} \ar[dl]^{\Id_{\B^{\ash}} \circ \epsilon_{\A^{\ash}}}\\ & \B^{\ash} &}.
$$
We therefore obtain the commutativity of diagrams (i) and (ii) and two surjections $p_{\A^{\ash}} : \A^{\ash} \vee^{\lambda} \B^{\ash} \rightarrow \A$ and $p_{\B^{\ash}} : \A^{\ash} \vee^{\lambda} \B^{\ash} \rightarrow \B$. Remembering the fact that the decomposition maps $\Delta_{\A^{\ash}}$, $\Delta_{\B^{\ash}}$ and $\Delta_{\A^{\ash} \vee^{\lambda} \B^{\ash}}$ are all the cofree decomposition map, we get the commutativity of the following diagrams:
$$
\mathrm{(a)} \hspace{1.5cm} \xymatrix@C=40pt@R=20pt{& \A^{\ash} \vee^{\lambda} \B^{\ash} \ar[d]^{\Delta_{\A^{\ash} \vee^{\lambda} \B^{\ash}}} \ar[dl]_{i_{2}} \ar[dr]^{i_{1}} & \\ \A^{\ash} \circ \B^{\ash} & \left(\A^{\ash} \vee^{\lambda} \B^{\ash}\right) \circ \left(\A^{\ash} \vee^{\lambda} \B^{\ash}\right) \ar[l]^{\hspace{-1cm} p_{\A^{\ash}} \circ p_{\B^{\ash}}} \ar[r]_{\hspace{1cm} p_{\B^{\ash}} \circ p_{\A^{\ash}}} & \B^{\ash} \circ \A^{\ash}}
$$
and
$$
\mathrm{(b)} \hspace{0.5cm} \xymatrix@R=16pt{\A^{\ash} \vee^{\lambda} \B^{\ash} \ar[d]_{\Delta_{\A^{\ash} \vee^{\lambda} \B^{\ash}}} \ar[r]^{i_{2}} & \A^{\ash} \circ \B^{\ash} \ar[r]^{\Id_{\A^{\ash}} \circ \Delta_{\B^{\ash}}} & \A^{\ash} \circ \B^{\ash} \circ \B^{\ash}\\ \left(\A^{\ash} \vee^{\lambda} \B^{\ash}\right) \circ \left(\A^{\ash} \vee^{\lambda} \B^{\ash}\right) \ar[r] & \left(\A^{\ash} \vee^{\lambda} \B^{\ash}\right) \circ \B^{\ash} \ar[ur]_{i_{2} \circ \Id_{\B^{\ash}}} &}
$$
and
$$
\mathrm{(c)} \hspace{0.5cm} \xymatrix@R=16pt{\A^{\ash} \vee^{\lambda} \B^{\ash} \ar[d]_{\Delta_{\A^{\ash} \vee^{\lambda} \B^{\ash}}} \ar[r]^{i_{1}} & \B^{\ash} \circ \A^{\ash} \ar[r]^{\Delta_{\B^{\ash}} \circ \Id_{\A^{\ash}}} & \B^{\ash} \circ \B^{\ash} \circ \A^{\ash}\\ \left(\A^{\ash} \vee^{\lambda} \B^{\ash}\right) \circ \left(\A^{\ash} \vee^{\lambda} \B^{\ash}\right) \ar[r] & \B^{\ash} \circ \left(\A^{\ash} \vee^{\lambda} \B^{\ash}\right) \ar[ur]_{\Id_{\B^{\ash}} \circ i_{1}} &}.
$$

Finally, the commutativity of the diagram (I) is a consequence of (a), (b) and (c), and of the coassociativity of $\Delta_{\A^{\ash} \vee^{\lambda} \B^{\ash}}$:
$$
\xymatrix@R=16pt@C=16pt{\B^{\ash} \A^{\ash} \ar[rr]^{i_{1}^{-1}} \ar[dd] && \A^{\ash} \vee^{\lambda} \B^{\ash} \ar[rr]^{i_{2}} \ar[dl] \ar[dr] && \A^{\ash} \B^{\ash} \ar[dd]\\
& \left(\A^{\ash} \vee^{\lambda} \B^{\ash}\right) \left(\A^{\ash} \vee^{\lambda} \B^{\ash}\right) \ar[d] && \left(\A^{\ash} \vee^{\lambda} \B^{\ash}\right) \left(\A^{\ash} \vee^{\lambda} \B^{\ash}\right) \ar[d] &\\
\B^{\ash} \B^{\ash} \A^{\ash} \ar[r]_{\Id_{\B^{\ash}} i_{1}^{-1}} & \B^{\ash} \left(\A^{\ash} \vee^{\lambda} \B^{\ash}\right) \ar[r]_{\Id_{\B^{\ash}} i_{2}} & \B^{\ash} \A^{\ash} \B^{\ash} \ar[r]_{i_{1}^{-1} \Id_{\B^{\ash}}} & \left(\A^{\ash} \vee^{\lambda} \B^{\ash}\right) \B^{\ash} \ar[r]_{i_{2} \Id_{\B^{\ash}}} & \A^{\ash} \B^{\ash} \B^{\ash}.}
$$
The case of the diagram (II) is similar.
$\cqfd$
\end{pf}

\begin{rem} \leavevmode
\begin{itemize}
\item We denote by $\A \vee_{\lambda} \B$ the operad $\P (V\oplus W,\, R \oplus D_{\lambda} \oplus S)$ so that the previous proposition provides a way to compute the decomposition map on $(\A \vee_{\lambda} \B)^{\ash} = \A^{\ash} \vee^{\lambda} \B^{\ash}$ in terms of $\Lambda^{c}$ and the decomposition maps on $\A^{\ash}$ and on $\B^{\ash}$ since
$$
\Delta_{\A^{\ash} \vee^{\lambda} \B^{\ash}} = (i_{1}^{-1} \circ i_{1}^{-1}) \cdot (\Id_{\B^{\ash}} \circ \Lambda^{c} \circ \Id_{\A^{\ash}}) \cdot (\Delta_{\B^{\ash}} \circ \Delta_{\A^{\ash}}) \circ i_{1}^{-1}.
$$
\item The diagrams (I) and (II) give a way to compute the map $\Lambda^{c}$ knowing the map $\lambda$.
\end{itemize}
\end{rem}

\subsection{The Diamond Lemma for distributive laws}\label{remdiamond}

The Diamond Lemma, Theorem 8.6.5 of \cite{LodayVallette}, provides an effective way of proving that $i_{1}$ is an isomorphism. Similarly as in the book \cite{LodayVallette}, the operads $\A$ and $\B$ and the cooperads $\A^{\ash}$ and $\B^{\ash}$ are weight graded by the opposite of the number of generators in $V$ and in $W$. Therefore, the $\Sb$-modules $\A \circ \B$, $\A \vee_{\lambda} \B$, $\B^{\ash} \circ \A^{\ash}$ and $\A^{\ash} \vee^{\lambda} \B^{\ash}$ are also weight graded by the opposite of the number of generators in $V$ minus the number of generators in $W$. The injective map $i_{1} : \A^{\ash} \vee^{\lambda} \B^{\ash} \rightarrow \B^{\ash} \circ \A^{\ash}$, resp. the surjective map $p : \A \circ \B \rightarrow \A \vee_{\lambda} \B$ given in \cite{LodayVallette}, preserves the weight grading and is always surjective, resp. injective, in weight $0$, $-1$ and $-2$. The following theorem is a slightly different version of Theorem 8.6.5 of \cite{LodayVallette} and says that if $p$ is an isomorphism in weight $-3$, the map $i_{1}$ is an isomorphism in any weight.

\begin{thm}\label{diamond}
Let $\A = \P(V,\, R)$ and $\B = \P(W,\, S)$ be two Koszul operads endowed with a rewriting rule $\lambda : W \circ_{(1)} V \rightarrow V \circ_{(1)} W$ such that the restriction of $p : \A \circ \B \twoheadrightarrow \A \vee_{\lambda} \B$ on $(\A \circ \B)^{(-3)}$ is injective. In this case, the morphisms $p$ and $i_{1} : \A^{\ash} \vee^{\lambda} \B^{\ash} \hookrightarrow \B^{\ash} \circ \A^{\ash}$ are isomorphisms, the map $\lambda$ induces a distributive law and the induced operad $(\A \circ \B,\, \gamma_{\Lambda})$ is Koszul, with Koszul dual cooperad $(\B^{\ash} \circ \A^{\ash},\, \Delta_{\Lambda^{c}})$.
\end{thm}

\begin{pf}
The only point to prove which is not in Theorem 8.6.5 of \cite{LodayVallette} is that $i_1$ is an isomorphism. In this case, we can conclude by Proposition \ref{decomposition map}. By Theorem 8.6.5 of \cite{LodayVallette}, we know that $p : \A \circ \B \cong \A \vee_0 \B \twoheadrightarrow \A \vee_{\lambda} \B$ is an isomorphism. It follows that the extension to the bar constructions $\Bm p : \Bm (\A \vee_0 \B) \rightarrow \Bm (\A \vee_{\lambda} \B)$ is also an isomorphism. On the bar construction, we consider a homological degree, called \emph{syzygy degree}, given by the weight degree minus the cohomological degree. The map $\Bm p$ is not dg but it commutes with the differential in syzygy degree $0$ up to a boundary given by $d_{\mathrm{bar}}(\lambda (w \otimes v))$ for $w \otimes v \in W \circ_{(1)} V$ and zero otherwise. By means of the fact that the bar constructions are zero in syzygy degree $1$, we get that $\Bm p$ descend to the syzygy degree $0$ homology group to give an isomorphism
$$
\Hm_0 \Bm p : \Hm_0 \Bm (\A \vee_0 \B) \cong (\A \vee_0 \B)^{\ash} \cong \A^{\ash} \vee^0 \B^{\ash} \cong \B^{\ash} \circ \A^{\ash} \rightarrow \Hm_0 \Bm (\A \vee_{\lambda} \B) \cong \A^{\ash} \vee^{\lambda} \B^{\ash}
$$
inverse to $i_1$.
$\cqfd$
\end{pf}

\begin{rem}
It is also enough to prove that $i_{1}$ is surjective in weight $-3$ to get the theorem. The proof is the same as the one in \cite{LodayVallette}, where we replace operads by cooperads, $\A \circ \B$ by $\B^{\ash} \circ \A^{\ash}$, $p$ by $i_{1}$, $\Bm$ by $\Omega$, the syzygy degree in \emph{Step 1} by the number of inversions and vice versa in \emph{Step 2}.
\end{rem}

\end{appendix}

% To include all references
%\nocite{*}
\bibliographystyle{alpha}
\bibliography{bib}

\def\cprime{$'$} \def\cprime{$'$} \def\cprime{$'$} \def\cprime{$'$}
  \def\cprime{$'$} \def\cprime{$'$}
\begin{thebibliography}{{Mne}08a}

\bibitem[BR73]{BernshteinRozenfeld}
I.~N. Bern{\v{s}}te{\u\i}n and B.~I. Rosenfel{\cprime}d.
\newblock Homogeneous spaces of infinite-dimensional {L}ie algebras and the
  characteristic classes of foliations.
\newblock {\em Uspehi Mat. Nauk}, 28(4(172)):103--138, 1973.

\bibitem[CFT02]{CattaneoFelderTomassini}
A.~S. {Cattaneo}, G.~{Felder}, and L.~{Tomassini}.
\newblock {From local to global deformation quantization of Poisson manifolds.}
\newblock {\em {Duke Math. J.}}, 115(2):329--352, 2002.

\bibitem[GK71]{GelfandKazhdan}
I.~M. Gel{\cprime}fand and D.~A. Ka{\v{z}}dan.
\newblock Certain questions of differential geometry and the computation of the
  cohomologies of the {L}ie algebras of vector fields.
\newblock {\em Dokl. Akad. Nauk SSSR}, 200:269--272, 1971.

\bibitem[{Gro}59]{Grothendieck}
A.~{Grothendieck}.
\newblock {G\'eom\'etrie formelle et g\'eom\'etrie alg\'ebrique. (Formal
  geometry and algebraic geometry).}
\newblock {Sem. Bourbaki 11 (1958/59), No.182, 28 p. (1959).}, 1959.

\bibitem[Hin97]{Hinich97}
V.~Hinich.
\newblock Homological algebra of homotopy algebras.
\newblock {\em Comm. Algebra}, 25(10):3291--3323, 1997.

\bibitem[{Hin}03]{Hinich03}
V.~{Hinich}.
\newblock {Erratum to ''Homological algebra of homotopy algebras''}.
\newblock {\em ArXiv Mathematics e-prints}, 2003.

\bibitem[HM12]{HirshMilles}
J.~{Hirsh} and J.~{Mill\`es}.
\newblock {Curved Koszul duality theory.}
\newblock {\em {Math. Ann.}}, 354(4):1465--1520, 2012.

\bibitem[Kad83]{Kadeishvili}
T.~Kadeishvili.
\newblock The algebraic structure in the homology of an ${A}( \infty
  )$-algebra.
\newblock {\em Soobshch. Akad. Nauk Gruzin. SSR}, no. 2(108):249--252, 1983.

\bibitem[Khu04]{Khudaverdian04}
H.~M. Khudaverdian.
\newblock Semidensities on odd symplectic supermanifolds.
\newblock {\em Comm. Math. Phys.}, 247(2):353--390, 2004.

\bibitem[{Kon}03]{Kontsevich03}
M.~{Kontsevich}.
\newblock {Deformation quantization of Poisson manifolds.}
\newblock {\em {Lett. Math. Phys.}}, 66(3):157--216, 2003.

\bibitem[LV12]{LodayVallette}
J.-L. Loday and B.~Vallette.
\newblock {\em Algebraic operads}, volume 346 of {\em Grundlehren der
  Mathematischen Wissenschaften}.
\newblock Springer, Heidelberg, 2012.

\bibitem[{Mer}04]{Merkulov04}
S.A. {Merkulov}.
\newblock {Operads, deformation theory and $F$-manifolds.}
\newblock In {\em {Frobenius manifolds. Quantum cohomology and singularities.
  Proceedings of the workshop, Bonn, Germany, July 8--19, 2002}}, pages
  213--251. 2004.

\bibitem[Mer05]{Merkulov05}
S.~A. Merkulov.
\newblock Nijenhuis infinity and contractible differential graded manifolds.
\newblock {\em Compos. Math.}, 141(5):1238--1254, 2005.

\bibitem[{Mer}06]{Merkulov06}
S.A. {Merkulov}.
\newblock {PROP profile of Poisson geometry.}
\newblock {\em {Commun. Math. Phys.}}, 262(1):117--135, 2006.

\bibitem[Mer08]{Merkulov08}
S.~A. Merkulov.
\newblock Lectures on {PROP}s, {P}oisson geometry and deformation quantization.
\newblock In {\em Poisson geometry in mathematics and physics}, volume 450 of
  {\em Contemp. Math.}, pages 223--257. Amer. Math. Soc., Providence, RI, 2008.

\bibitem[Mer10]{Merkulov10}
S.~A. Merkulov.
\newblock Wheeled {P}ro(p)file of {B}atalin-{V}ilkovisky formalism.
\newblock {\em Comm. Math. Phys.}, 295(3):585--638, 2010.

\bibitem[{Mne}08a]{Mnev08b}
P.~{Mnev}.
\newblock {Discrete BF theory}.
\newblock {\em ArXiv e-prints}, 2008.

\bibitem[Mne08b]{Mnev08a}
P.~Mnev.
\newblock On the simplicial {$BF$}-theory.
\newblock {\em Dokl. Akad. Nauk}, 418(3):308--312, 2008.

\bibitem[MV09]{MerkulovVallette2}
S.~{Merkulov} and B.~{Vallette}.
\newblock {Deformation theory of representations of prop(erad)s. II.}
\newblock {\em {J. Reine Angew. Math.}}, 636:123--174, 2009.

\bibitem[Nij55]{Nijenhuis}
A.~Nijenhuis.
\newblock Jacobi-type identities for bilinear differential concomitants of
  certain tensor fields. {I}, {II}.
\newblock {\em Nederl. Akad. Wetensch. Proc. Ser. A. {\bf 58} = Indag. Math.},
  17:390--397, 398--403, 1955.

\bibitem[NN57]{NewlanderNirenberg}
A.~Newlander and L.~Nirenberg.
\newblock Complex analytic coordinates in almost complex manifolds.
\newblock {\em Ann. of Math. (2)}, 65, 1957.

\bibitem[Sch93]{Schwarz93}
Albert Schwarz.
\newblock Geometry of {B}atalin-{V}ilkovisky quantization.
\newblock {\em Comm. Math. Phys.}, 155(2):249--260, 1993.

\bibitem[Str09]{Strohmayer}
H.~Strohmayer.
\newblock Operad profiles of {N}ijenhuis structures.
\newblock {\em Differential Geom. Appl.}, 27(6):780--792, 2009.

\bibitem[{Str}10]{Strohmayer10}
H.~{Strohmayer}.
\newblock {Prop profile of bi-Hamiltonian structures.}
\newblock {\em {J. Noncommut. Geom.}}, 4(2):189--235, 2010.

\end{thebibliography}
\ \\
\noindent
\textsc{Joan Millès, Universit\'e de Paul Sabatier, Institut de mathématiques de Toulouse, 118, route de Narbonne, 31062 Toulouse Cedex 9, France}\\
E-mail : \texttt{joan.milles@math.univ-toulouse.fr}\\
\textsc{URL :} \texttt{http://www.math.univ-toulouse.fr/$\sim$jmilles}\\

\end{document}